\documentclass[letterpaper, 12 pt]{article}
\usepackage{amsthm, amsmath, amssymb, amsfonts, graphicx}

\newcommand{\setleftmargin}[1]{
        \addtolength{\textwidth}{\oddsidemargin}
        \addtolength{\textwidth}{1in}
        \addtolength{\textwidth}{-#1}
        \setlength{\oddsidemargin}{-1in}
        \addtolength{\oddsidemargin}{#1}
        \setlength{\evensidemargin}{\oddsidemargin}
}
\newcommand{\setrightmargin}[1]{
        \setlength{\textwidth}{8.5in}
        \addtolength{\textwidth}{-\oddsidemargin}
        \addtolength{\textwidth}{-1in}
        \addtolength{\textwidth}{-#1}
}
\setleftmargin{1 in}
\setrightmargin{1 in}

\newcommand{\ZZ}{\mathbb{Z}}

\newcommand{\RR}{\mathbb{R}}

\newcommand{\R}{\mathcal{R}}

\newcommand{\MM}{\mathcal{M}}
\newcommand{\kk}{\kappa}

\newcommand{\equivalent}{\Longleftrightarrow}

\DeclareMathOperator{\Trop}{Trop}

\DeclareMathOperator{\In}{In}

\DeclareMathOperator{\RowSpan}{RowSpan}
\DeclareMathOperator{\diag}{diag}

\DeclareMathOperator{\Hull}{Hull}
\DeclareMathOperator{\link}{link}

\DeclareMathOperator{\Spec}{Spec}

\newcommand{\DD}{\mathcal{D}}
\newcommand{\EE}{\mathcal{E}}
\newcommand{\oDD}{\mathring{\DD}}

\newcommand{\remove}{\backslash \backslash}

\newcommand{\scap}{\underset{\textrm{stable}}{\cap}}

\newtheorem{conj}{Conjecture}[section]
\newtheorem{thrm}[conj]{Theorem}
\newtheorem{theorem}[conj]{Theorem}
\newtheorem{prop}[conj]{Proposition}
\newtheorem{Lemma}[conj]{Lemma}
\newtheorem{lemma}[conj]{Lemma}

\newtheorem{cor}[conj]{Corollary}

\newtheorem*{f-vector}{The $f$-Vector Conjecture}
\newtheorem*{f-vector-sp}{The $f$-Vector Conjecture, Continued}

\theoremstyle{definition}
\newtheorem{problem}[conj]{Problem}
\newtheorem{remark}[conj]{Remark}

\begin{document}

\title{Tropical Linear Spaces}
\author{David E Speyer \\ \texttt{speyer@math.berkeley.edu}}
\maketitle

\begin{abstract}
We define tropical analogues of the notions of linear space and Pl\"ucker coordinate and study their combinatorics. We introduce tropical analogues of intersection and dualization and define a tropical linear space built by repeated dualization and transverse intersection to be constructible. Our main result that all constructible tropical linear spaces have the same $f$-vector and are ``series-parallel''. We conjecture that this $f$-vector is maximal for all tropical linear spaces with equality precisely for the series-parallel tropical linear spaces. We present many partial results towards this conjecture.

In addition we describe the relation of tropical linear spaces to linear spaces defined over power series fields and give many examples and counter-examples illustrating aspects of this relationship. We describe a family of particularly nice series-parallel linear spaces, which we term tree spaces, that realize the conjectured maximal $f$-vector and are constructed in a manner similar to the cyclic polytopes.
\end{abstract}

\section{Introduction and Summary}

The last five years have shown an explosion of interest in the field of
tropical geometry. In tropical geometry, multiplication is replaced by
addition and addition by minimum. See \cite{TropMath} or \cite{FirstSteps} for an introduction.
A polynomial
is said to vanish tropically if, among the values of its terms, the
minimum is not unique. For example, $xy+x+y$ tropically vanishes if $x=y
\leq x+y$ or $x=x+y \leq y$ or $y=x+y \leq x$. One motivation for this
definition is that this is the union of the three sets $x+y=\min(x,y)$,
$x=\min(x+y,y)$ and $y=\min(x,x+y)$ which can be obtained by replacing
multiplication by addition and addition by minimum in the various
equations $-xy=x+y$, $-x=xy+y$ and $-y=xy+x$. Another justification for
this process is that it replaces exact polynomial relations between
variables by approximate relations between the negatives of their logs.
(We have $-\log(xy)=-\log x - \log y$ and $-\log(x+y)=\min(-\log x, - \log
y)+O(1)$ for $x$ and $y>0$.) A third, most precise, justification is that
if we study power series solutions $X_i(t)$ to a polynomial $f$ then
tropicalizing this polynomial gives the necessary and, it turns out (see
\cite{SpeySturm}, theorem 2.1) also sufficient relations between the
vanishing orders of the
$X_i(t)$. We will discuss this perspective further in section
\ref{nonarch}.

This purpose of this paper is to study the tropical analogues of the
notion of linear spaces and Pl\"ucker coordinates. One of the
pleasures of this theory is that it can be set up in a purely
combinatorial manner without reference to general tropical
theory. However, doing so requires defining slightly more inclusive
notions of tropical linear space than if the theory were exactly set
up to match the classical geometry. This is very similar to the
situation of matroid theory, where the definition of a matroid is
simple and every operation that can be performed with a vector
arrangement has a corresponding operation for a general matroid. On
the other hand, characterizing those matroids that come from vector
arrangements -- the so-called realizable matroids -- is very
hard. Similarly, we will have a notion of realizable tropical linear
space which will be associated to a linear space over a
nonarchimedean field and most of our motivation will come from
operations with linear spaces.  Our approach, however, will be to
first give purely combinatorial definitions and only then discuss the
relation to geometry.

Tropical linear spaces will be polyhedral complexes. Thus, we can ask
about all of the invariants associated to polyhedral complexes. Ardila
and Klivans \cite{ArdKliv} have studied the link of a vertex in a tropical
linear space and have shown it to be homeomorphic to the
chain complex of the lattice of flats of a certain matroid and thus,
(\cite{White2}, sect. 7.6), homotopic to a wedge of spheres. The
polytopes occuring in tropical linear spaces are Minkowski summands of
permutahedra
(proposition \ref{permut}). The topology of tropical linear spaces is quite
simple: they
are contractible (theorem \ref{contract}). So we understand the local
combinatorics
and the local
and global toplogy of tropical linear spaces.

The global combinatorics of tropical linear spaces, on the other hand, is
quite intriguing. There is a great deal of theoretical and some
experimental evidence for the following:
\begin{f-vector}
The number of $i$-dimensional faces of a tropical $d$-plane in
$n$-space which become bounded after being mapped to
$\RR^n/(1,\ldots,1)$ is at most $\binom{n-2i}{d-i}
\binom{n-i-1}{i-1}$.
\end{f-vector}

\begin{remark}
It would follow from the $f$-vector conjecture that the number of total $i$-dimensional faces of a tropical $d$-plane in $n$-space, without any boundedness condition, is at most $\binom{n-i-1}{d-i} \binom{2n-d-1}{i-1}$. See proposition \ref{blah}.
\end{remark}

The following is a summary of the rest of the paper.

In section \ref{basic}, we will define the basic concepts mentioned in
this introduction and prove their essential properties. We will also
introduce an operation called dualization, which is analogous to the
classical operation of taking the dual, or orthogonal complement, of a
linear space.

In the section \ref{intersect}, we will introduce an operation of
transverse
intersection that produces new tropical linear spaces from old. We define
a tropical linear space to be \emph{constructible} if it can be built from
tropical hyperplanes (equivalently, from points) by successive dualization
and transverse intersection. Much of the rest of the paper will be devoted
toward proving:

\begin{theorem} \label{cons-f}
Every constructible space achieves the $f$-vector of the $f$-vector
conjecture.
\end{theorem}

In section \ref{nonarch}, we pause to explain how tropical linear spaces
are
related to linear spaces over power series fields, just as matroids are
related to linear spaces over ordinary fields. These section includes
many counter-examples to show that tropical linear spaces can fail to be
realizable in almost every concievable way.

In section \ref{sp} we will define a notion of ``series-parallel
tropical linear space'', analagous to the notion of series-parallel
matroid. We also describe a way of working with series-parallel matroids
in
terms of two-colored trees that will be important in the future. One of
our themes will be that series-parallel tropical linear spaces are
the best tropical linear spaces. In particular we conjecture
\begin{f-vector-sp}
Equality in the $f$-vector conjecture is achieved precisely by
series-parallel linear spaces.
\end{f-vector-sp}
We prove
\begin{theorem} \label{cons-sp}
Every constructible space is series-parallel.
\end{theorem}

One unfortunate consequence of theorem \ref{cons-sp} is that it
is quite hard to find a general method for writing down tropical
linear spaces that are not series-parallel! More precisely, it is
not so hard to produce degenerate limits of series-parallel linear spaces
which are not themselves series-parallel. Writing down a large number
of tropical linear spaces which are not limits of series-parallel
tropical linear spaces, however, is a challenge -- which is one of the
reasons that experimentally testing the $f$-vector conjecture is
tricky.

In section \ref{cases}, we will prove the $f$-vector conjecture in the
cases $i=1$ and $d=n/2$. We then return to proving our main results. We
prove theorem \ref{cons-sp} in section \ref{results} and also prove
many lemmas that will be used in the proof of theorem \ref{cons-f}. We
also prove theorem \ref{cons-f} in section \ref{results}.

Finally, in section \ref{trees}, we will introduce a notion of tree
linear space which achieves the bounds in the $f$-vector conjecture
and has very explicit combinatorics.  Suggestively, the construction
of tree linear spaces will be reminiscent of the construction of
cyclic polytopes.

\section{Basic Definitions and Results} \label{basic}

For $f=\sum_{e \in E} f_e X^e$ a nonzero polynomial, $E \subset
\ZZ_{\geq 0}^n$, we denote by $\Trop(f)$ the set of all points $(x_1,
\ldots, x_n)$ in $\RR^n$ such that, if we form the collection of
numbers $\sum_{i=1}^n e_i x_i$ for $e$ ranging over $E$, then the
minimum of this collection is not unique. Throughout the paper, we will use $X$, $P$ and occasionally $Q$ to refer to variables in classical algebro-geometric coordinates and $x$, $p$ and $q$ for the corresponding tropical variables.

Consider a collection of variables $p_{i_1 \ldots i_d}$ indexed by the $d$
element subsets of $[n]:=\{ 1, \ldots, n \}$. We will say that $p$ is a
\emph{tropical Pl\"ucker vector} if
$$(p_I) \in \Trop (P_{Sij} P_{Skl} - P_{Sik} P_{Sjl} + P_{Sil} P_{Sjk})$$
for every $S \in \binom{[n]}{d-2}$ and every $i$, $j$, $k$ and $l$
distinct members of $[n] \setminus S$.

\begin{remark}
This definition is equivalent to saying that $p_I$ is a
valuated matroid with values in the semiring $(\RR, +, \min)$. See
\cite{Dress1} as well as the more general \cite{Dress2} and \cite{Dress3}.
\end{remark}

Let $\Delta(d,n)$ denote the $(d,n)$-hypersimplex, defined as the
convex hull of the points $e_{i_1} + \cdots + e_{i_d} \in \RR^n$ where
$\{ i_1, \ldots, i_d \}$ runs over $\binom{[n]}{d}$. We abbreviate
$e_{i_1}+\cdots+e_{i_d}$ by $e_{i_1 \ldots i_d}$.  Consider a real-valued
function $\{ i_1, \ldots, i_d \} \mapsto p_{i_1 \ldots i_d}$ on the
vertices of $\Delta(d.n)$. We define a
polyhedral subdivision $\DD_p$ of $\Delta(d,n)$ as follows: consider the
points $(e_{i_1}+\cdots+e_{i_d}, p_{i_1 \cdots i_d}) \in \Delta(d,n)
\times \RR$ and take their convex hull. Take the lower faces (those
whose outward normal vector have last component negative) and project
them back down to $\Delta(d,n)$, this gives us the subdivision $\DD_p$.
We will often omit the subscript $p$ when it is clear from context. A
subdivision which is obtained in this manner is called \emph{regular},
see, for example, \cite{Zieg} definition 5.3.

Let $P$ be a subpolytope of $\Delta(d,n)$. We say that $P$ is matroidal if
the vertices of $P$, considered as elements of $\binom{[n]}{d}$ are the
bases of a matroid $M$. In this case, we write $P=P_M$. Our references for
matroid terminology and theory are Neil White's anthologies \cite{White1}
and
\cite{White2}.

\begin{prop}
The following are  equivalent:

\begin{enumerate}
\item $p_{i_1 \ldots i_d}$ are tropical Pl\"ucker coordinates

\item The one skeleta of $\DD$ and $\Delta(d,n)$ are the same.

\item Every face of $\DD$ is matroidal.
\end{enumerate}
\end{prop}

\begin{proof}
$(1) \implies (2)$. Every edge of $\DD$ joins two vertices of
$\Delta(d,n)$. If $e$ is an edge of $\DD$ connecting $e_I$ and $e_J$, we
define the \emph{length} of $e$, $\ell(e)$, to be $|I \setminus J|=|J
\setminus I|$. Our claim is that every edge of $e$ has length $1$. We
prove by induction on $\ell \geq 2$ that $\DD$ has no edge of length
$\ell$.

For the base case, if $e$ is an edge with $\ell(e)=2$ then $e=(e_{Sij},
e_{Skl})$ for some $S \in \binom{[n]}{d-2}$. The six vertices $e_{Sij}$,
$e_{Sik}$, $e_{Sil}$, $e_{Sjk}$, $e_{Sjl}$, $e_{Skl}$ form the vertices of
an octahedron $O$ with $e_{Sij}$ and $e_{Skl}$ opposite vertices. One can
check that the condition $(p_I) \in \Trop (P_{Sij} P_{Skl} - P_{Sik}
P_{Sjl} + P_{Sil} P_{Sjk})$ implies that $O$ is either a face of $\DD$ or
subdivided in $\DD$ into two square pyramids (in one of three possible
ways). In any of these cases, $e$ is not an edge of $\DD$.

Now consider the case where $\ell>2$. Suppose (for contradiction) that $e$
is an edge of $\DD$. Let $e=(e_{S T}, e_{S T'})$, $T \cap T' =\emptyset$,
$|T|=|T'|=\ell$. Let $F$ be the face of $\Delta(d,n)$ consisting of all
vertices $e_I$ with $S \subset I \subset S \cup T \cup T'$. Then $e$ must
belong to some
two dimensional face of $\DD$ contained in $F$; call this two dimensional
face $G$.

Let $\gamma$ be the path from $e_{ST}$ to $e_{S T'}$ that goes around
$G$ the other way from $e$. No two vertices of $F$ are more than
distance $\ell$ apart, so the edges of $\gamma$ have lengths less than
or equal to $\ell$. If $\gamma$ contained an edge $(e_{SU}, e_{SU'})$
of length $\ell$ then its midpoint $(e_{SU}+e_{SU'})/2$ would also be
the midpoint of $e$ contradicting the convexity of $G$. Thus, every
edge of $\gamma$ has length less than $\ell$ and by induction must have
length
$1$. So all the edges of $\gamma$ are in the direction $e_i-e_j$ for
some $i$ and $j \in [n]$. These vectors must span a two dimensional
space. This means either that there are $\{ i_1, i_2, j_1, j_2 \}
\subset [n]$ such that all edges of $\gamma$ are parallel to some
$e_{i_r}-e_{j_s}$ or there are $\{ i_1, i_2, i_3 \} \subset [n]$ such
that every edge of $\gamma$ is parallel to $e_{i_r}-e_{i_s}$. In either
case, $e$ has length at most $2$, but that returns us to our ground
case.

$(2) \implies (3)$: Let $P$ be any polytope in $\DD$. By assumption, all
of the edges of $P$ are edges of $\Delta(d,n)$. It is a theorem of
Gelfand, Goresky, MacPherson and Serganova (\cite{GGMS}, theorem 4.1) that
this implies that $P$ is
matroidal.
Since the proof is short, we include it. Let $e_I$ and $e_J$ be vertices
of $P$ with $j \in J \setminus I$. We must prove that there is a vertex of
$P$ of the form $e_{I \cup \{ j \} \setminus \{ b \}}$ for some $b \in I$.

Define a linear functional $\phi : \Delta(d,n) \to \RR$ by $\phi(x_1,
\ldots x_n)=\sum_{i \in I} x_i + d x_j$. Then $Q := P \cap \{ x: \phi(x)
\geq d \}$ is a convex polytope and hence connected. $Q$ contains the
vertices $e_I$ and $e_J$ and there is a path from $e_I$ to $e_J$ along
edges of $P$ which lie in $Q$. Let the first step of this path go from
$e_I$ to $e_{I \cup \{ a \} \setminus \{ b \}}$. If $a \neq j$, $\phi(e_{I
\cup \{ a \} \setminus \{ b \}})=d-1$ and $e_{I \cup \{ a \} \setminus \{
b \}}$ does not lie in $Q$. So instead we have $a=j$ and we are done.

$(3) \implies (1)$: It is easy to check that, if $(1)$ is false, $\DD$ has
a one dimensional face of the form $\Hull(e_{Sij}, e_{Skl})$, with $i$,
$j$, $k$ and $l$ distinct. This is not matroidal.
\end{proof}

Now, suppose that $(p_I) \in \RR^{\binom{[n]}{d}}$ obeys the tropical
Pl\"ucker relations. Define $L(p) \subset \RR^n$ by
$$\bigcap_{1 \leq j_1 < \cdots < j_{d+1} \leq n} \Trop(\sum_{r=1}^{d+1}
(-1)^r P_{j_1 \cdots \hat{j_r} \cdots j_{d+1}} X_{j_r}).$$
We term any $L$ which arises in this manner a \emph{$d$-dimensional
tropical linear space in $n$-space}. We often omit the dependance on $p$
when it is clear from context.

$L(p)$ is essentially the same set as  the tight span of $p$ defined in
\cite{TreeOfLife}. There are two differences: (1) Dress's sign conventions
are opposite to ours and (2) $L(p)$ is invariant under translation by
$(1,\ldots,1)$, Dress chooses a particular representative within each
orbit for this translation.

While the above definition makes the connection to ordinary linear spaces
most clear, for practically every purpose it is better to work with the
alternate characterization which we now give. For any $w \in \RR^n$,
define $\DD_w$ to be the subset of the vertices of $\Delta(d,n)$ at which
$p_{i_1 \ldots i_d} - \sum_{j=1}^d w_{i_j}$ is minimal. This is, by
definition, a face of $\DD$ and thus is $P_{M_w}$ for a matroid $M_w$.

\begin{prop}
$w \in L(p)$ if and only if $M_w$ is loop-free.
\end{prop}

Recall that a matroid is called loop-free if every element of the matroid
appears in at least one basis. There is a geometric way of recognizing
when $M$ is loop-free from the polytope $P_M$: $M$ is loop-free if and only if $P_M$
is not contained in any of the $n$ facets of $\Delta(d,n)$ of the type
$x_i =0$ for $1 \leq i \leq n$. In particular, note that if $P_M$ meets
the interior of $\Delta(d,n)$, $M$ is necessarily loop-free.

\begin{proof}
By replacing $p_{i_1 \ldots i_d}$ by $p_{i_1 \ldots i_d} - \sum_{j=1}^d
w_{i_j}-(\textrm{constant})$, we may assume without loss of generality
that
$w=0$ and $\min p_I=0$. Then
$M_w=M_0$ is simply the matroid whose basis correspond to the $I$ for
which $p_I=0$. First, we assume that $M_0$ has a loop $j$ and prove that
$0 \not \in L(p)$. Let $(i_1, \ldots, i_d)$ be a basis of $M_0$, clearly
$j \not\in (i_1, \ldots, i_d)$. Then $p_{i_1 \ldots i_d}=0$ but $p_{j i_1
\ldots \hat{i_{r}} \ldots i_d} >0$ for $1 \leq r \leq d$. Taking $(j_1,
\ldots, j_{d+1})=(j, i_1, \ldots, i_d)$ we see that $0 \not \in \Trop(
\sum_{r=1}^{d+1} (-1)^r p_{j_1 \ldots \hat{j_r} \ldots j_{d+1}} x_{j_r})$.

The converse is more interesting. Fix $J=\{ j_1, \ldots, j_{d+1} \}$, our
aim will be to prove $0 \in \Trop( \sum_{r=1}^{d+1} (-1)^r p_{j_1 \ldots
\hat{j_r} \ldots j_{d+1}} x_{j_r})$. Let $e_J=\sum e_{j_r}$ and, for any
$s \in \RR$, set $M_s=M_{-s e_J}$. Note that, for $s$ large enough, all of
the bases of $M_s$ will be subsets of $J$. It is equivalent to show that,
for such an $s$, $M_s$ has at least two bases or, equivalently, that for
any $j \in J$ and any such $s$, $j$ is not a loop of $M_s$. By hypothesis,
$j$ is not a loop of $M_0$, so it is enough to show that if $j \in J$ is a
loop of $M_s$ for some $s$ then it is a loop of $M_{s'}$ for all $s' < s$.

As $s$ varies, $M_s$ changes at a finite number of values of $s$,
call them $s_1 < s_2 < \cdots < s_k$. Suppose that $s_i < s < s_{i+1}$,
then $P_{M_s}$ is a face of both $P_{M_{s_i}}$ and of $P_{M_{s_{i+1}}}$.
The bases of $M_s$ are precisely the bases of $M_{s_i}$ which have the
largest possible number of elements in common with $J$. Similarly, the
bases of $M_s$ are  precisely the bases of $M_{s_{i+1}}$ which have the
smallest possible number of elements in common with $J$. In other words,
$$M_s = M_{s_i}|_J \oplus M_{s_i}/J = M_{s_{i+1}}|_{[n] \setminus J}
\oplus
M_{s_{i+1}}/([n] \setminus J).$$

From the displayed equation, it follows that if $j
\in J$ is a loop in $M_{s_{i+1}}$ then it is a loop in
$M_s$.
Similarly, if $j$ is a loop in $M_s$ then it is a loop in
$M_{s_i}$. Concatenating deductions of this sort, we see that, as
promised, if $j$ is a loop of $M_s$ then it is a loop of $M_{s'}$ for all
$s' < s$.
\end{proof}

As $\Delta(d,n)$ is contained in the hyperplane $x_1 + \cdots + x_n=d$
we see $L$ is invariant under translation by $(1,1,\ldots,1)$. We will
abuse notation by saying a face of $L$ is bounded if its image in
$\RR^n / (1,\ldots, 1)$ is bounded and calling a face a vertex if its
image in $\RR^n/(1,\ldots,1)$ is zero dimensional. However, when we
refer to the dimensions of faces of $L$ we will always be speaking of
$L$ itself, without taking the quotient by $(1,\ldots,1)$.

\begin{figure}
\centerline{\includegraphics{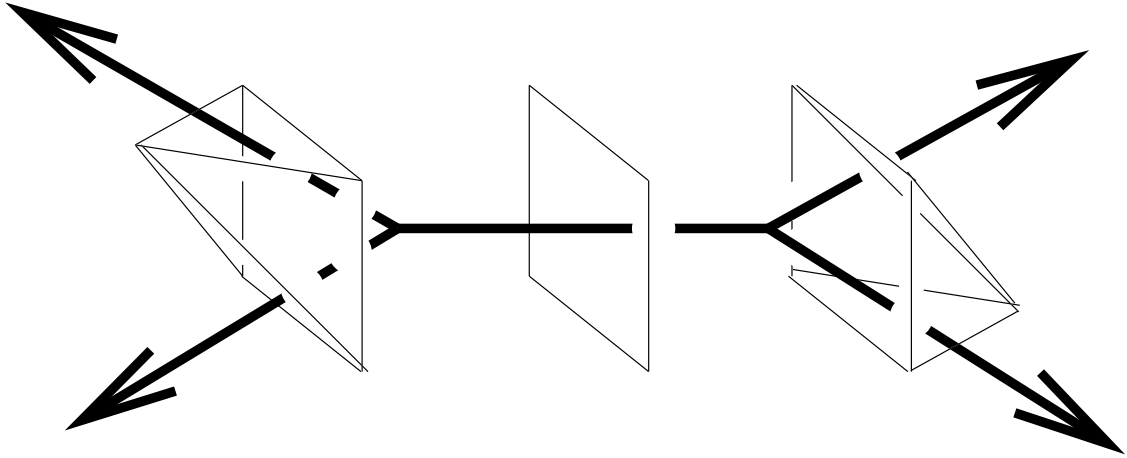}}
\caption{A Subdivision of $\Delta(2,4)$ and the Corresponding $L/(1,1,1,1)$}\label{octahedron}
\end{figure}

We see that $L$ is a subcomplex of $\DD^{\vee}$, where
$\DD^{\vee}$ is defined to be the polyhedral subdivision of $\RR^n$
where $w$ and $w'$ lie in the same face if $M_w=M_{w'}$. In figure \ref{octahedron}, we show $\Delta(2,4)$, which is an octahedron, subdivided into two square pyramids and  draw  the dual $L/(1,1,1,1)$ in bold. We write
$M^{\vee}$ for the face of $\DD^{\vee}$ dual to $P_M$. If
$v=M^{\vee}/(1,\ldots,1)$ is a vertex of $L/(1,\ldots,1)$ we see that the
link of $L/(1,\ldots,1)$ at $v$ is the subcomplex of the normal fan to
$P_M$ consisting of the normals to the loop-free facets. This fan is
studied in \cite{ArdKliv}, we summarize the main results of that paper in
the
following proposition:

\begin{prop}
$\link_{M^{\vee}}(L)$ is homeomorphic to the chain complex of the lattice
of flats of $L$. If $N^{\vee}$ is a face of $L$ containing $v$ then
$N=\bigoplus F_{i+1}/F_i$ for some flag of flats $\emptyset = F_0 \subset
F_1 \subset \cdots \subset F_r=[n]$ of $M$.
\end{prop}

When $M^{\vee}/(1,\ldots,1)$ is positive dimensional, $P_M$ turns out to
be a product of $P_N$'s.

\begin{prop} \label{Sum}
Let $M$ be a matroid with $M=\bigoplus M_i$ the decomposition of $M$ into
connected parts. Then $P_M = \prod P_{M_i}$.
\end{prop}

\begin{proof}
By the definition of direct sum, the bases of $M$ are all combinations of
bases of the $M_i$.
\end{proof}

We can now describe the faces of $L$. Recall that the braid arrangement
$\mathcal{B}_r$ is the fan in $\RR^{r}$ whose facets are the cones of the
type $x_{\sigma(1)} \leq \ldots \leq x_{\sigma(r)}$ where $\sigma$ ranges
over the symmetric group $S_r$.

\begin{prop} \label{permut}
If $M^{\vee}$ is an $r$-dimensional face of $L$ then the normal fan to
$M^{\vee}$ is a subfan of a coarsening of $\mathcal{B}_r$. When $M^{\vee}$
is bounded, $M^{\vee}$ is a Minkowski summand of the $r^{\textrm{th}}$
permutahedron.
\end{prop}

\begin{proof}
The normal fan to $M^{\vee}$ is the same as the link of $P_M$ in $\DD$.
Let $P_{\tilde{M}}$ be a face of $\DD$ with $P_M$ a face of
$P_{\tilde{M}}$. We must show that the local cone of $P_{\tilde{M}}$ at
$P_M$ is a union of cones of $\mathcal{B}_r$. But  the local cone of
$P_{\tilde{M}}$ at $P_M$ is precisely the negative of the cone of the
normal fan to $P_{\tilde{M}}$ dual to $P_M$. From the description in
\cite{ArdKliv}, we see that the normal fan to $P_{\tilde{M}}$ is a
coarsening of $\mathcal{B}_n$ and that the cones dual to codimension $i$
faces are unions of cones of $\mathcal{B}_r$.

If $M^{\vee}$ is bounded then its normal fan is complete and thus, by the
last paragraph, a coarsening of $\mathcal{B}_r$. $\mathcal{B}_r$ is the
normal fan to the $r^{\textrm{th}}$ permutahedron. If $P$ and $Q$ are two
polytopes such that the normal fan to $Q$ refines that of $P$ then $P$ is
a Minkowski summand of $tQ$ for $t$ sufficiently large.
\end{proof}

Not every coarsening of $\mathcal{B}_r$ is the normal fan to a polytope.
For example, take a cube, divide one of its faces into two triangles and
cone from the center. This is a coarsening of $\mathcal{B}_4$ which is not the normal fan of any polytope. I do not
know whether there are further obstructions to a fan showing up as the
normal fan to a face of a tropical linear space.

\begin{problem}
Let $P$ be a polytope whose normal fan is a coarsening of $\mathcal{B}_r$.
Can $P$ always occur as a face of a tropical linear space?
\end{problem}

For ordinary
linear spaces, the linear space determines the Pl\"ucker coordinates up to
scaling. For tropical linear spaces, a similar result holds.

\begin{prop}
Let $L$ be a tropical linear space with tropical Pl\"ucker coordinates
$p_I$. Then $L$ determines the $p_I$ up to addition of a common scalar.
\end{prop}

\begin{proof}
Let $S$ be a $d-1$ element subset of $[n]$ and let $F$ denote the face of
$\Delta(d,n)$ whose vertices are those $I$ containing $S$. Since $F$ is an
$(n-d)$-simplex, it can not be subdivided in $\DD$ so $F$ is a face of
$\DD$. $F$ corresponds to a loop-free matroid $\Phi$ so it is dual to a
face
$\Phi^{\vee}$ of $L$. $\Phi^{\vee}$ is of the form $(\RR_{\geq 0}^S +
\RR(1,\ldots,1))+q$ for $q$ a vector satisfying $q_i=p_{Si}$ for $i
\not\in S$. No other face of $L$ is of the form $(\RR_{\geq 0}^S \oplus
\RR(1,\ldots,1))+q$ and $q$ is unique up to translation by $(1,\ldots,1)$.
Thus, $L$ determines $p_{Si}-p_{Sj}$. As this is true for all $S$, $i$ and
$j$, $L$ determines the $p_I$ up to adding the same constant to all of
them.
\end{proof}

\begin{theorem} \label{contract}
$L$ is a pure $d$-dimensional contractible polyhedral complex.
\end{theorem}

\begin{proof}
To show that $L$ is pure $d$-dimensional, it is enough to show that the
link of any vertex $v$ of $L/(1,\ldots,1)$ is a pure
$(d-2)$-dimensional complex. Let $v+(1,\ldots,1)$ be dual to $P_M$ for
some matroid
$M$. By the results of \cite{ArdKliv}, the link of $v$ has a subdivision
isomorphic to the order complex of the lattice of flats of $M$. This
lattice is graded of length $d$, so its order complex is pure of dimension
$d-2$.

We now show contractibility. It is easy to see that each tropical
hyperplane is tropically convex in
the sense of \cite{DevStur}. Thus, their intersection is tropically
convex and hence contractible.
\end{proof}

We now study the bounded faces of $L$.

\begin{prop}
The following are equivalent:

\begin{enumerate}

\item $w$ lies in a bounded face of $L$

\item $P_{M_w}$ is an interior face of $\DD$.

\item $M_w$ is loop-free and co-loop-free.

\end{enumerate}
\end{prop}

\begin{proof}
$(2) \equivalent (3)$: $P_{M_w}$ is interior if and only if it is not contained
in any facet of $\Delta(d,n)$. The facets of $\Delta(d,n)$ are of two
types: $x_i=1$ and $x_i=0$. $P_{M_w}$ lies in the former type if and only if $i$
is a co-loop of $M_w$; $P_{M_w}$ lies in the latter type if and only if $i$ is a
loop of $M_w$.

$(1) \equivalent (2)$: If $P$ is not an interior face of $\DD$, the
corresponding dual face is not bounded. Conversely, if $P$ is internal,
the corresponding dual face of $\DD^{\vee}$ is bounded, so we just must
check that this dual face is in $L(p)$ at all. For this, we must check
that $M_w$ is loop-free. But we saw in the previous paragraph that this
follows from $P_{M_w}$ being internal.
\end{proof}

At this point, we can prove our earlier remark.

\begin{prop} \label{blah}
If the $f$-vector conjecture is true then a tropical $d$-plane in $n$-space has no more than $\binom{n-i-1}{d-i} \binom{2n-d-1}{i-1}$ faces of dimension $i$.
\end{prop}

\begin{proof}
Let $\DD$ be a matroidal decomposition of $\Delta(d,n)$, we must bound the number of codimension $i-1$ faces $P_M$ of $\DD$ such that $M$ has no loops. For such a $P_M$, let $\Lambda \subset [n]$ be the co-loops of $M$ and let $\Delta$ be the face of $\Delta(d,n)$ whose vertices are the sets $I\in \binom{[n]}{d}$ which contain $\Delta$. The claim that $M$ has no loops is equivalent to saying that $P_M$ is an interior face of the decomposition of $\Delta$. So we may count the total number of $P_M$ by summing over $\Delta$. 

There are $\binom{n}{j}$ faces $\Delta$ for which $|\Lambda|=j$; each such $\Delta$ is isomorphic to $\Delta(d-j,n-j)$. A subpolytope of such a $\Delta$ has codimension $i-1$ in $\Delta(d,n)$ if and only if it has co-dimension $i-j-1$ in $\Delta$. So, assuming the $f$-vector conjecture and writing the bound iun the $f$-vector conjecture as $\binom{n-i-1}{d-i} \binom{n-d-1}{i-1}$, we have the following bound for the number of co-dimension $i-1$ loop-free faces of $\DD$:

\begin{eqnarray*}
\sum_j \binom{n}{j} \binom{(n-j)-(i-j)-1}{(d-j)-(i-j)} \binom{(n-j)-(d-j)-1}{(i-j)-1} &=& \\
 \binom{n-i-1}{d-i} \sum_j \binom{n}{j} \binom{n-d-1}{i-j-1} &=&  \binom{n-i-1}{d-i} \binom{2n-d-1}{i-1}.
\end{eqnarray*}
The last sum is evaluated by the identity $\sum_{j} \binom{p}{j} \binom{q}{r-j}=\binom{p+q}{r}$. 
\end{proof}

Suppose that $p: \binom{[n]}{d} \to \RR$ obey the tropical Pl\"ucker
relations. Define $p^{\perp} : \binom{[n]}{n-d} \to \RR$ by
$p^{\perp}_I=p_{[n] \setminus I}$. Then the subdivision $\DD^{\perp}$ of
$\Delta(n-d,n)$ induced by $p^{\perp}$ is simply the image of $\DD$ under
$(x_1, \ldots, x_n) \mapsto (1-x_1, \ldots, 1-x_n)$. This map replaces
each matroidal polytope with the polytope of the dual matroid; in
particular, $p^{\perp}$ also obeys the tropical Pl\"ucker relations.

We denote $L(p^{\perp})$ by $L^{\perp}$. We can easily check that $L
\mapsto L^{\perp}$ is an inclusion reversing bijection from the set of
tropical linear spaces in $\RR^n$ to itself.

\begin{prop}
The bounded part of $L^{\perp}$ is negative that of $L$. If $L$ is a
tropical linear space of dimension $d$ in $n$ space then the
bounded part of $L$ is
at most $\min(d,n-d)$-dimensional.
\end{prop}

\begin{proof}
The first sentence follows because the property of being loop- and
co-loop-free is self dual. $L$ is $d$-dimensional and $L^{\perp}$ is
$(n-d)$-dimensional; as the bounded part of $L$ is isomorphic to a
subcomplex of both it is at most $\min(d,n-d)$-dimensional.
\end{proof}

\section{Stable Intersections} \label{intersect}

If $L$ and $L'$ are two ordinary linear spaces of dimensions $d$ and $d'$
in $n$ space which meet transversely then $L \cap L'$ is a $(d+d'-n)$-dimensional space whose Pl\"ucker coordinates are given by
$$P_J(L \cap L')=\sum_{\substack{I \cap I'=J \\ |I|=d \\ |I'|=d'}} \pm
P_I(L) P_{I'}(L').$$
Note that the summation conditions on $I$ and $I'$ guarantee that $I \cup
I'=[n]$.

We will discuss a tropical version of this formula. If $L$ and $L'$ are
tropical linear spaces of dimension $d$ and $d'$, define a real valued
function $p(L \scap L')$ on $\binom{[n]}{d+d'-n}$ by
$$p_J(L \scap L')=\min_{\substack{I \cap I'=J \\ |I|=d \\ |I'|=d'}}
(p_I(L)+p_{I'}(L')).$$

At the moment, $p_J(L \scap L')$ is just defined as a formal symbol,
there is not yet any tropical linear space called $L \scap L'$. We now
prove that there is such a space.

\begin{prop}
$p(L \scap L')$ is a tropical Pl\"ucker vector.
\end{prop}

\begin{proof}
Set $P$ to be the polyhedron of points above the lower convex hull of the
points $(p_{i_1 \ldots i_d}(L), e_{i_1} + \cdots + e_{i_d}) \in \RR \times
\Delta(d,n)$ and define $P' \subset \RR \times \Delta(d',n)$ similarly.
Include $\Delta(d+d'-n,n)$ into the Minkowski sum $\Delta(d,n) +
\Delta(d',n)$ by $\iota : e \to e+(1,\ldots, 1)$. Set $Q=(P + P') \cap
(\RR
\times \iota(\Delta(d+d'-n,n)))$. Then projecting the bottom faces of $Q$
back
down to $\Delta(d+d'-n,n)$ yields the regular subdivision of
$\Delta(d+d'-n,n)$ induced by $p(L \scap L')$.

So we see that every face of $\DD_{p(L \scap L')}$ is of the form
$$R:=\iota( \Delta(d+d'-n,n) ) \cap \left( P_M + P_{M'}
\right).$$
The vertices of this face are the points $\iota(e_J) \in
\iota(\Delta(d+d'-n,n))$ where $J=I \cap I'$ for $I$ and $I'$ bases
of $M$ and $M'$ with $I \cup I'=[n]$.

There is an operation $M \wedge M'$, called matroid intersection
(\cite{White1}, section 7.6), which takes two matroids of ranks $d$
and $d'$ and produces a third matroid. The spanning sets of $M \wedge
M'$ are precisely the sets of the form $U \cap U'$ where $U$ and $U'$
span $M$ and $M'$ respectively. $M \wedge M'$ always has rank at least
$d+d'-n$. If the rank of $M \wedge M'$ is larger than
$d+d'-n$ then there are no pairs of bases $I$ and $I'$ of $M$ and $M'$
with $I \cup I'=[n]$.  In this case $R$ is empty and does not
contribute a face to $\DD_{p(L \scap L')}$. Alternatively, the rank of
$M \wedge M'$ is $d+d'-n$. Then we see that $R=P_{M \wedge M'}$.
\end{proof}

We have the immediate corollary

\begin{cor}
Every face of $\DD_{p(L \scap L')}$ is of the form $P_{M \wedge M'}$ for
$P_M$ and $P_{M'}$ faces of $\DD_{p(L)}$ and $\DD_{p(L')}$.
\end{cor}

We thus see that there is a well defined tropical linear space $L \scap
L'$. Our aim now is to justify the notation $L \scap L'$ by showing that
$L \scap L'$ is the ``stable'' intersection of $L$ and $L'$. We first need
a definition and a combinatorial lemma:

Let $M$ and $M'$ be two loop-free matroids on $[n]$ with connected
decomposition $M=\bigoplus_{i \in I} M_i$, $M'=\bigoplus_{i \in I'}
M'_{i}$. Let  $S_i$ and $S'_i$ be the ground sets of $M_i$ and $M'_i$. Let
$\Gamma(M,M')$ be the bipartite graph (possibly with multiple edges) whose
vertex set is $I \sqcup I'$ and which has an edge between $I$ and $I'$
for each element of $S_i \cap S'_{i'}$. We say that $M$ and $M'$ are
\emph{transverse} if $\Gamma(M,M')$ is a forest without multiple edges.

\begin{Lemma}\label{BasisExists}
Suppose that $M$ and $M'$ are transverse and let $\Gamma=\Gamma(M,M')$. Then $M \wedge M'$ is a loop-free matroid whose connected components are in bijection with the components of $\Gamma$. $M \wedge M'$ is the matroid formed by repeatedly taking parallel connections along the edges of $\Gamma$. $M \wedge M'$ has rank $d+d'-n$. 
\end{Lemma}

The phrase ``repeatedly taking parallel connections'' deserves a proper definition. What we mean is the following: let $F$ be a forest, let $E$ be a finite set and let each edge of $F$ be labelled with an element of $E$, no two edges receiveing the same label. Suppose that, for each vertex $v \in F$, we are given a subset $S(v )$ of  $E$. We require that $S(v) \cap S(w) = \{ e \}$ if $e$ is the edge $(v,w)$, that $S(v) \cap S(w)=\emptyset$ if there is no edge joining $v$ and $w$ and that $\bigcup_{v \in F} S(v)=E$. For every $v \in F$, we place the structure of a connected loop-free matroid $M(v)$ on $S(v)$. We then define a matroid $N$ on the ground set $E$ by the following procedure:

\begin{enumerate}
\item If $F$ has no edges, let $N=\bigoplus_{v \in F} M(v)$.

\item Otherwise, let $v$ and $w$ be vertices of $F$ joined by an edge $e$. 

\item Form a new forest $F'$ where $e$ is contracted to a single vertex $u$. Set $S'(u)=S(v) \cup S(w)$ and let $M'(u)$ be the parallel connection of $M(v)$ and $M(w)$ along $e$. For all $t \in F$ other than $u$, set $S'(t)=S(t)$ and $M'(t)=M(t)$.

\item Return to step (1) with the new forest $F'$ and the new $S'(\cdot)$ and $M'(\cdot)$.

\end{enumerate}

The claim is that $N$ does not depend on the order of the contractions and, that when we use $F=\Gamma$, $S(i)=S_i$ or $S'_i$, for $i \in I $ or $I'$ and $M(i)=M_i$ or $M'_i$ for $i \in I$ or $I'$ then $N=M \wedge M'$.

\begin{proof}
Let $T$, $S(\bullet)$ and $M(\bullet)$ be as above and let $N(F,S(\cdot),M(\cdot))$ be the a matroid constructed by the above procedure. We claim that $U \subset E$ is a basis of $N(F,S(\bullet),M(\bullet))$ if and only if there exist bases $B(v)$ of $M(v)$ such that, for every $e \in E$, precisely one of the following four conditions holds:
\begin{enumerate}
\item $e \in U$, $e=(v,w)$ is an edge of $F$ and $e \in B(v) \cap B(w)$.
\item $e \in U$, $e \in S(v)$ is not an edge of $F$ and $e \in B(v)$. 
\item $e \not \in U$,  $e=(v,w)$ is an edge of $F$ and $e$ lies in precisely one of $B(v)$ and $B(w)$.
\item  $e \not \in U$, $e \in S(v)$ is not an edge of $F$ and $e \not \in B(v)$\end{enumerate}

This is the definition of direct sum when $F$ has no edges and it is easy to check that this description is unchanged by contracting an edge. Note that this description is independent of the order of contractions, so we see that $N(T, S(\bullet), M(\bullet))$ is independent of the order of contraction. 

Applying the above description in the case where $F=\Gamma(M,M')$ \emph{etc.}, we see that the bases of $N(\Gamma, M_{\bullet}, S_{\bullet})$ are precisely the sets of the form $B \cap B'$ where $B$ and $B'$ are bases of $M$ and $M'$ with $B \cup B'=[n]$. These are precisely the bases of $M \wedge M'$ providing that some such $(B, B')$ exist. But the iterative construction above clearly does produce a matroid and every matroid has at least one basis, so such $B$ and $B'$ do exist. $|B \cap B'|=|B| \cup |B'| - | B \cup B'|$ so $N$ has rank $d+d'-n$. Fianlly, observe that parallel connections preserves loop-freeness, so $N$ is loop-free.
\end{proof}

\begin{prop}
Let $0 \leq d, d' \leq n$ with $d+d' \geq n$. Suppose that $L$ and $L'$
are tropical linear spaces in $n$-space of dimension $d$ and $d'$. Then $L \scap L'
\subseteq L \cap L'$. If $L$ and $L'$ meet transversely then $L \scap L'=L
\cap L'$.
\end{prop}

\begin{proof}
Let $\DD$, $\DD'$ and $\EE$ be the subdivisions of $\Delta(d,n)$,
$\Delta(d',n)$ and $\Delta(d+d'-n,n)$ corresponding to $L$, $L'$ and $L \scap L'$. Let
$\DD^{\vee}$, $(\DD')^{\vee}$ and $\EE^{\vee}$ denote the dual
subdivisions of $\RR^n$.

Consider a particular $w \in \RR^n$, we need to show that $w \in L
\scap L'$ implies that $w \in L$ and $w \in L'$ and that the reverse
holds if $L$ and $L'$ meet transversely. Let $M^{\vee}$, $(M')^{\vee}$
and $N^{\vee}$ be the faces of $\DD^{\vee}$, $(\DD')^{\vee}$ and
$\EE^{\vee}$ respectively containing $w$ and let $P_M$, $P_{M'}$ and $P_N$
be
the respective dual faces of $\DD$, $\DD'$ and $\EE$.  Then $P_N=(P_M
+ P_{M'}) \cap \iota(\Delta(d+d'-n,n))$. We must show that, if
$M$ or $M'$ has a loop then $N$ has a loop and that, if $L$ and $L'$
meet transversely and $N$ has a loop then $M$ or $M'$ does as well.

First, suppose that $M$ contains a loop $e$. Then, for every $x \in
P_M$, $x_e=0$. For every $x \in P_{M'}$, $0 \leq x_e \leq 1$. Thus,
for $x \in P_M + P_{M'}$, $0 \leq x_e \leq 1$ and thus $x_e=0$ on
$P_N$. So $e$ is a loop of $N$.

Now assume that $L$ and $L'$ meet transversely. Let $M = \bigoplus_{i \in
I} M_i$ and $M'=\bigoplus_{i \in I'} M'_i$ be the components of $M$ and
$M'$ and let $S_i$ and $S'_i$ be the ground sets of $M_i$ and $M'_i$
respectively.

\begin{Lemma} \label{isforest}
$M$ and $M'$ are transverse, \emph{i.e.}, $\Gamma(M,M')$ is a forest and has no multiple edges.
\end{Lemma}

\begin{proof}
Let $|I|=e$ and $|I'|=e'$, so $e$ and $e'$ are the dimensions of
$M^{\vee}$ and $(M')^{\vee}$. The affine linear space spanned by
$M^{\vee}$ is cut out by the equations $x_i-x_j= \textrm{constant}$ when
edges $i$ and $j$ have the same endpoint in $I$, and similarly for
$(M')^{\vee}$. Thus, $M^{\vee} \cap (M')^{\vee}$ spans the affine
linear space cut out by the equations $x_i - x_j=\textrm{constant}$
whenever edges $i$ and $j$ are in the same component of $\Gamma$. We
assumed that $L$ and $L'$ meet transversely, so the dimension of the last
linear space must be $e+e'-n$ and $\Gamma$ must have $e+e'-n$ connected components. A graph with $e+e'$ vertices, $n$ edges and
$e+e'-n$ connected components must be a forest without multiple edges.
\end{proof}

We now know by lemma \ref{BasisExists} that $M \wedge M'$ is loop-free. 
\end{proof}

The next theorem explains the motivation for the notation $L \scap L'$ and
the geometric meaning of $L \scap L'$ when $L$ and $L'$ are not
trasnverse.

\begin{theorem}
Let $L$ and $L'$ be tropical linear spaces. Then, for a generic $v \in
\RR^n$, $L$ and $L'+v$ meet transversely. $L \scap L'=\lim_{v \to 0} (L
\cap (L'+ v))$.
\end{theorem}

\begin{proof}
Let $F$ and $F'$ be faces of $L$ and $L'$. Suppose that the affine linear
spaces spanned by $F$ and $F'$ together fail to span $\RR^n$. Then, for
$v$ outside a proper subspace of $\RR^n$, $F \cap (F' + v)=\emptyset$.
Thus, for a generically chosen $v$, every such $F$ and $F'+v$ fail to
meet. The only faces of $L$ and $L'+v$ that do meet, then, meet
transversely.

We have $p_I(L'+v)=p_I(L')+\sum_{i \in I} v_{i}$ so the tropical Pl\"ucker
coordinates of $L'+v$ vary continuously with $v$. The tropical Pl\"ucker
coordinates of $L \scap (L'+v)$ similarly vary continuously with $v$. When
$v$ is chosen generically, $L$ and $L'+v$ meet transversely so $L \cap
(L'+v) = L \scap (L' + v)$. Take limits of both sides as $v \to 0$ and use
the proceeding continuity arguments to conclude that $\lim_{v \to 0} (L
\cap (L'+v))=L \scap L'$.
\end{proof}

\section{Linear Spaces over Power Series Fields} \label{nonarch}

Let $\kk$ be an algebraically closed field and let $K$ be an algebraic
closure of $\kk((t))$. With some care as to our language, all of our
results can be generalized to a general algebraically closed field
complete with respect to a nontrivial archimedean valuation but we
will stick to this case for simplicity. There is a unique valuation
$v:K^* \to \RR$ extending $v(a t^{\alpha} + \cdots)=\alpha$.  Let $R$
be the ring $\{0 \} \cup v^{-1}(\RR_{\geq 0})$ and $\MM$ the ideal $\{
0 \} \cup v^{-1}(\RR_{>0})$ so $R/\MM=\kk$.  If $X$ is a subvariety of
$K^n$, we define $\Trop X$ to be the closure in $\RR^n$ of $v(X \cap
(K^*)^n)$. We also, when it is more convenient, write $\Trop I$ where
$I$ is the ideal of $X$.

Let $f \in K[x_1, \ldots, x_n]$ be a nonzero polynomial and let $w \in
\RR^n$. We define a polynomial $\In_w f \in \kk[x_1, \ldots, x_n]$ as
follows: expand $f(t^{w_1} x_1, \ldots, t^{x_n} x_n)$ as a power series
$\sum t^{\alpha_i} f_i(x_1, \ldots, x_n)$ where $f_i \in \kk[x_1,
\ldots, x_n] \setminus \{ 0 \}$ and $\alpha_i$ is increasing. Then $\In_w
f=f_0$.  If
$I \subset K[x_1, \ldots, x_n]$ is an ideal, we define $\In_w I$ to be
the ideal generated by $\In_w f$ for $f \in I$. Similarly, if $X
\subset K^n$ is a variety with ideal $I$, we denote $\Spec
(\kk[x]/\In_w I)$ as $\In_w X$.

If $w \in \Trop I$ and $f \in I \setminus \{ 0 \}$ then $\In_w f$ is clearly not a monomial. It turns
out that the converse is true as well.

\begin{thrm}
$\Trop I$ is the set of $w$
for which $\In_w I$ contains no monomial.
\end{thrm}

\begin{proof}
See \cite{SpeySturm}, theorem 2.1.
\end{proof}

We earlier defined $\Trop f$ for $f$ a polynomial and have just now given a definition of $\Trop (f)$ for $(f)$ a principal ideal. For $f \in \kk[x]$, these two definitions are consistent.

Now consider the case where $L$ is a linear space with Pl\"ucker
coordinates $P_{i_1 \ldots i_d}$. Let $p_{i_1 \ldots i_d}=v(P_{i_1 \ldots
i_d})$. The $P$'s obey the Pl\"ucker relations, so the $p$'s obey the
tropical Pl\"ucker relations. Our next proposition shows that the
combinatorially defined $L(p)$ truly does reflect the geometrically
defined $\Trop L$.

\begin{prop}
$\Trop L=L(p)$.
\end{prop}

\begin{proof}
Every point $(x_1, \ldots, x_n)$ in $L$ obeys $\sum (-1)^r P_{j_1 \ldots
\widehat{j}_r \ldots j_{d+1}} X_{j_r}=0$ so the $v(X_i)$ lie in
$\Trop(\sum_{r=1}^{d+1}
(-1)^r P_{j_1 \cdots \hat{j_r} \cdots j_{d+1}} X_{j_r})$. Thus, $\Trop L
\subseteq L(p)$.

For the converse direction, suppose $w \in L(p)$, we will prove that
$w \in \Trop L$. Without loss of generality, let $w=0$ and suppose
that $0 = \min_{I \in \binom{[n]}{d}} p_I$. Let $P_I(0) \in \kk$ be
the image of $P_I$ in $\kk=R/\MM$. Clearly, the $P_I(0)$ obey the
Pl\"ucker relations and hence correspond to an ideal $L(0) \subset
\kk^n$. It is easy to see that $L(0)=\In_w L$. Let $M$ be the set of
$I$ for which $p_I=0$, by assumption $M$ is the set of bases of a loop
free matroid. So $L(0)$ is not contained in any of the coordinate
planes of $\kk^n$ and thus the corresponding ideal contains no
monomial.
\end{proof}

It also turns out that $\scap$ does actually capture the behavior of
linear spaces under intersection.

\begin{prop}
Let $L$ and $L' \subset K^n$ be two linear spaces. Then there exists a
linear space $\tilde{L'} \subset K^n$ with $\Trop \tilde{L'}=\Trop L'$ and $\Trop L
\scap \Trop L' = \Trop L \scap \Trop \tilde{L'}=\Trop (L \cap
\tilde{L'})$. If $\Trop L$ and $\Trop L'$ meet transversely, then we
already have $\Trop L
\scap \Trop L'=\Trop (L \cap L')$ without having to choose an
$\tilde{L'}$.
\end{prop}

\begin{proof}
Let $d$ and $d'$ be the dimensions of $L$ and $L'$. Let $P_I$, $P'_I$ and $Q_I$ be the Pl\"ucker coordinates of $L$, $L'$ and
$L \cap L'$, so
$$Q_J=\sum_{\substack{I \cap I'=J \\ |I|=d \\ |I'|=d' }}
\pm P_I P'_{I'}.$$
If $q_I$ are the tropical Pl\"ucker coordinates of $\Trop L \scap \Trop
L'$, we have
$$q_J=\min_{\substack{ I \cap I'=J \\ |I|=d \\ |I'|=d' }}
v(P_I P'_{I'}).$$
We want to have $q_J=v(Q_J)$, so we will be successful if there is no
cancellation of leading terms in $\sum \pm P_I P'_{I'}$.

Let $(u_1, \ldots, u_n) \in (\kk^*)^n$. Let $\tilde{L'}=\diag(u_1,
\ldots, u_n) L'$. Letting $\tilde{P'_I}$ denote the Pl\"ucker
coordinates of $\tilde{L'}$, we have $\tilde{P'_I}= P'_I \prod_{i \in
I} u_i$. When the $u_i$ are chosen generically, there is no
cancellation of leading terms in the sum $\sum \pm P_I \tilde{P'_{I'}}
= \sum \pm P_I P'_I \prod_{i \in I'} u_i$.  For a generic $u$,
therefore, $\tilde{L'}$ has the desired property.

We now must prove the second claim, that if $\Trop L$ and $\Trop L'$ meet
transversely then we never have cancellation in this sum. We will do this
by showing that one term has a lower valuation than all of the others. Let
$\EE$ be the subdivision of $\Delta(d+d'-n,n)$ induced by $q$ and let $\DD$
and $\DD'$ be the subdivisions of $\Delta(d, n)$ and $\Delta(d',
n)$ corresponding to $\Trop L$ and $\Trop L'$. Let $P_N$ be a facet of
$\EE$
containing the vertex $J$ and let $w \in N^{\vee}$. Let $M^{\vee}$ and
$(M')^{\vee}$ be the faces of $\DD^{\vee}$ and $(\DD')^{\vee}$
containing $w$. Without loss of generality, we can take $w=0$ and $0=\min
p_I=\min p'_I$.

Let $\bigoplus_{i \in I} M_i$ and $\bigoplus_{i \in I'} M'_i$ be the
connected components of $M$ and $M'$. Write $\Gamma=\Gamma(M,M')$. As
in lemma \ref{isforest}, we see that $\Gamma$ is a forest. The fact
that $J$ is a basis of $ N$ indicates that there are $B$ and $B'$ with $B \cap
B'=J$, $B \cup B'=[n]$ such that, for each $M_i$, $B \cap M_i$ a basis
of $M_i$ and similarly for $B'$. If we show there is only one such
pair $(B, B')$, that will indicate that there is only one term in the
sum with valuation zero and hence no cancellation.

Suppose there were two such pairs, $(B_1, B'_1)$ and $(B_2, B'_2)$. For $e
\in [n]$ and $S \subset [n]$, define $[e \in S]$ to be $1$ if $e \in S$
and $0$ otherwise. Note that $[e \in B_s] + [e \in B'_s] = 1 + [e \in J]$.
Set
\begin{eqnarray*}
\alpha(e) &=& [e \in B_1] - [e \in B'_1] - [e \in B_2] + [e \in B'_2] \\
&=& 2 \left( [e \in B_1] - [e \in B_2] \right).
\end{eqnarray*}
we think of $\alpha$ as a cochain in $C^1(\Gamma)$. We have
\begin{eqnarray*}
\partial \alpha(M_i) &=& 2 \sum_{e \in M_i}  \left( [e \in B_1] - [e \in
B_2] \right) \\
&=& 2 ( |B_1 \cap M_i| - |B_2 \cap M_i|).
\end{eqnarray*}
Since $B_1 \cap M_i$ and $B_2 \cap M_i$ are both bases of $M_i$ they have
the same cardinality and $\partial \alpha(M_i)=0$ for every $M_i$.
Similarly, $\partial \alpha(M'_i)=0$. So $\alpha$ is a cochain. But
$\Gamma$ is acyclic, so $\alpha=0$. This in turn implies that  $[e \in
B_1] = [e \in B_2]$ for every $e \in [n]$ so $B_1=B_2$ and similarly
$B'_1=B'_2$.
\end{proof}

Having seen that every linear space over $K$ gives rise to a tropical
linear space, one might
ask whether every tropical linear space arises in this manner. The answer
to this question is
a dramatic ``no'' and we conclude this section by showing various manners
in which this can
fail. We will say that a tropical linear space is realizable if it arises
from a linear space
over $K$.

First, $\DD$ can contain polytopes corresponding to non-realizable
matroids; such a polytope
obviously can not occur in a realizable decomposition. For example, if $M$
is the non-Pappus
matroid, removing $P_M$ from $\Delta(3,9)$ leaves $8$ matroidal polytopes,
each of them a
cone on $\Delta_2 \times \Delta_5$. Cutting $\Delta(3,9)$ into the
non-Pappus matroid and the
eight other pieces gives a non-realizable matroidal subdivision. More
generally, any matroid
at all may appear as a piece of a matroidal subdivision:

\begin{prop}
If $M$ is a rank $d$ matroid on $n$ elements, the function $- \rho(\cdot)$
on
$\binom{[n]}{d}$ obeys the tropical Pl\"ucker relations, where $\rho(S)$
is the rank of the
flat spanned by $S$ for any $S \subset [n]$. $M$ is a face of the
corresponding subdivision.
\end{prop}

\begin{proof}
Let $S \in \binom{[n]}{d-2}$ and let $i$, $j$, $k$ and $l$ be distinct
elements of $[n]
\setminus S$. Write $\rho_{M/S}$ for the rank function of
$M/S$, we have
$\rho(Sij)=\rho_{M/S}(ij)+\rho(S)$. Thus, we only need to check that
$-\rho$ obeys the
tropical Pl\"ucker relations in the case where $d=2$ and $n=4$, which is straight
forward. $M$ is a face
of the corresponding subdivision because $-\rho$ is minimal on the
vertices of $M$.
\end{proof}

\begin{figure}
\centerline{\includegraphics{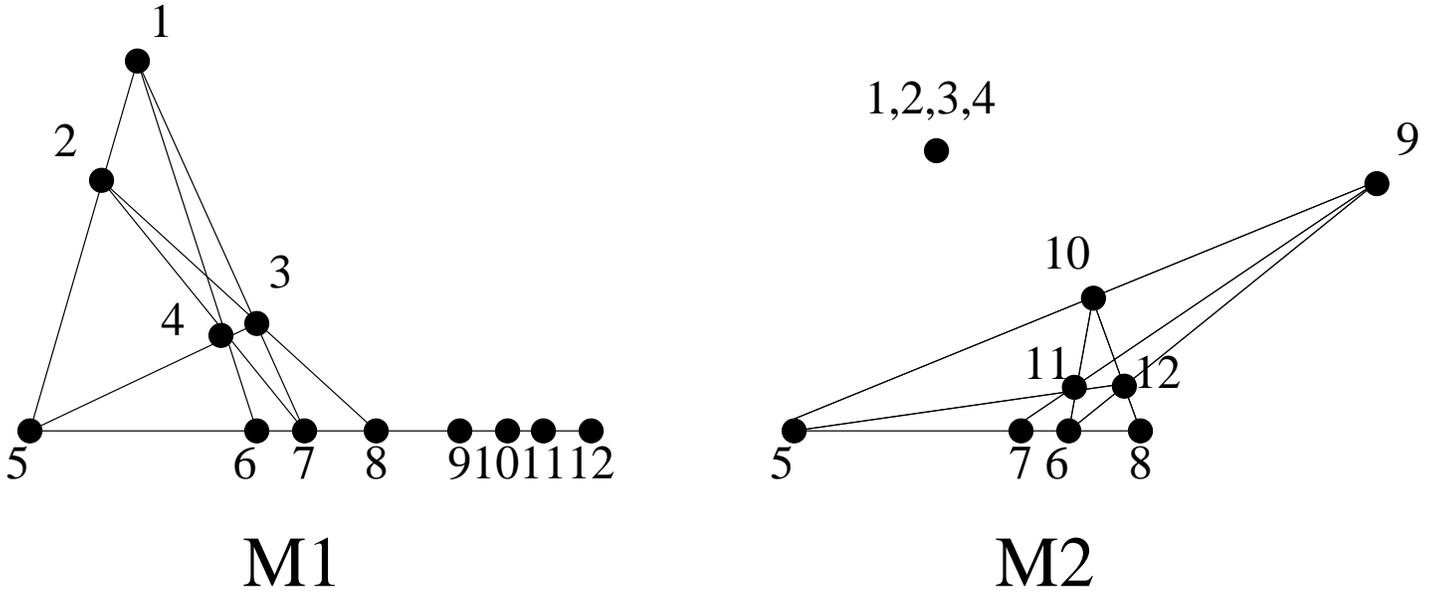}}
\caption{The Matroids $M_1$ and $M_2$} \label{CrossRatio2}
\end{figure}

Even if every face of $\DD$ is realizable, it is still possible for
the subdivision to be non-realizable for global reasons. Suppose $\kk$
does not have characteristic $2$ or $3$. Let $M_1$ and $M_2$ be the rank 3
matroids on 12 points corresponding to the plane geometries in figure
\ref{CrossRatio2}.  After removing $M_1$ and $M_2$ from
$\Delta(3,12)$, the remainder can be cut into cones on $\Delta_2
\times \Delta_8$. This gives a matroidal subdivision into realizable
matroids.  This subdivision is not realizable; if $p_{ijk} \in K$
where the Pl\"ucker coordinates of a linear space over $K$ realizing
this subdivision then the presence of $M_1$ implies that $p_{156}
p_{178} = p_{167}p_{158} (1+\textrm{higher\ order\ terms})$. The
presence of $M_2$ implies that $p_{156} p_{178} = p_{167} p_{158}
(-2+\textrm{higher\ order\ terms})$, a contradiction.

\begin{figure}
\centerline{\includegraphics{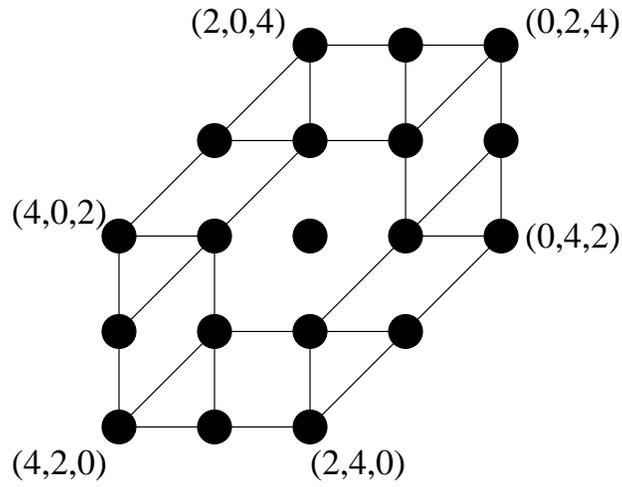}}
\caption{A non-regular subdivsion that induces a non-regular subdivision
of $\Delta(6,12)$.} \label{NonReg}
\end{figure}

It is possible to have a matroidal subdivision that is not regular and
hence certainly can not be realizable. Map $\Delta(6,12)$ to
$\ZZ^3/(1,1,1)$ by $(x_1, \ldots, x_{12}) \to (x_1+\ldots+x_4,
x_5+\ldots+x_8, x_9+\ldots+x_{12})$. Subdivide $\Delta(6,12)$
according to the preimage of figure \ref{NonReg}. The resulting
subdivision is not regular because figure \ref{NonReg} is not.

\begin{figure}
\centerline{\includegraphics{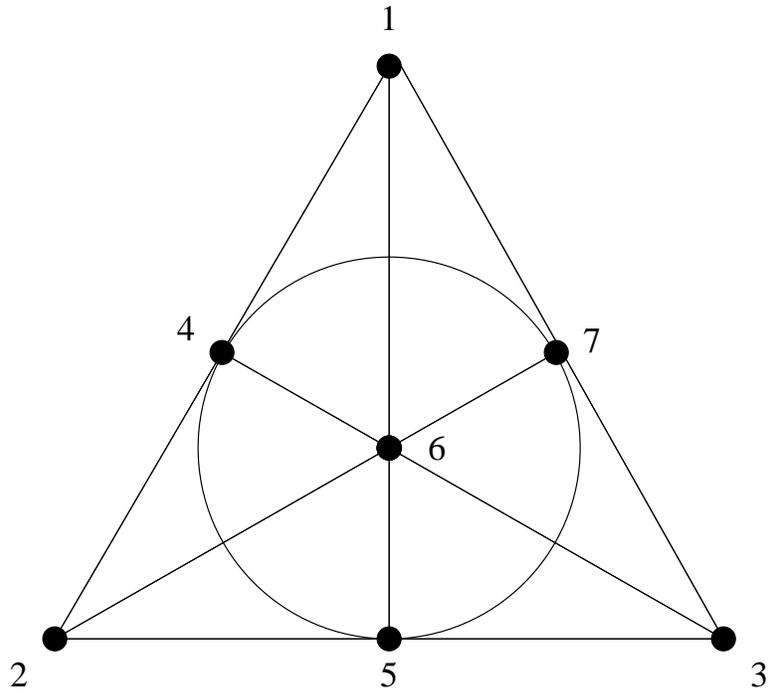}}
\caption{The Fano Plane} \label{Fano}
\end{figure}

It is possible for a matroidal subdivision to be induced by two
different $p$'s such that one of the $p$'s is realizable
and the other is not. We give an example when $K$ has
characteristic $2$, but there are characteristic zero examples. Let
$F$ denote the Fano plane (see figure \ref{Fano}), and let $\DD$ be
the subdivision of $\Delta(3,7)$ into $P_F$ and $7$ cones on $\Delta_2
\times \Delta_3$. Let $C \subset \binom{[7]}{3}$ be the triples of
dependent elements of $F$.  Let $p$ be a collection of tropical
Pl\"ucker coordinates inducing $\DD$.  We can assume without loss of
generality that $p$ is $0$ on $P_F$. Then $p$ induces $\DD$ whenever
$p_I>0$ for $I \in C$.

In characteristic $2$, the following relation holds among the Pl\"ucker
coordinates:
\begin{eqnarray*} & 0 \quad =  \quad
\,-\, P_{367} P_{567} \underline{P_{124}}
- P_{167} P_{467} \underline{P_{235}}
- P_{127} P_{567} \underline{P_{346}}  \\ & \quad
- P_{126} P_{367} \underline{P_{457}}
- P_{237} P_{467} \underline{P_{156}}
+ P_{134} P_{567} \underline{P_{267}}
+ P_{246} P_{567} \underline{P_{137}}
 \,\, + \, P_{136} \underline{P_{267}} \underline{P_{457}}.
\end{eqnarray*}
Here the underlined variables are those corresponding to the elements of
$C$. We see that,
with our asssumption that $p$ is zero on $P_F$ and in particular on
all vertices outside $C$, the valuations of the various terms are
$p_I$ as $I$ runs over $C$ and also $p_{267}+p_{457}$. As $p_{457}>0$,
this last can't be minimal. We see that the $p's$ must obey the
relation that the minimum of $p_I$ for $I \in C$ is not unique  in order
to be realizable. (In
face, this condition is sufficient for realizability as well.)

We end by discussing two naturally occuring examples of tropical
Pl\"ucker coordinates for which I do not know whether or not the
corresponding linear space is realizable.

\begin{problem}
Let $T$ be a tree with leaves
labelled by $[n]$ and a positive weight on each edge. For $I \subset [n]$,
define $[I]$ to be the minimal subtree of $T$ containing $I$ and set $p_I$
to be negative the sum of the weights of all edges in $[I]$. It is shown
in \cite{SpeyLior} that $p_I$ are tropical Pl\"ucker coordinates, are they
realizable?
\end{problem}

\begin{problem}
Let $\R=\RR((t))$, $\R$ is an ordered field where $a
t^{\alpha}+ \cdots$ is positive if and only if $a$ is. Thus, it makes
sense to define a $d \times d$ matrix $A$ with entries in $\R$ to be
positive definite if $vAv^T>0$ for all $v \in \R^d \setminus \{ 0
\}$. Let $X_1$, \ldots, $X_n$ be $d \times d$ positive definite
symmetric matrices with entries in $\R$. Set $P_{i_1 \ldots i_d}$
equal to the coefficient of $x_{i_1} \cdots x_{i_d}$ in $\det \sum x_i
X_i$. The $P_I$ do not obey the Pl\"ucker relations. Nonetheless, it
is shown in \cite{Horn} that $v(P_I)$ do obey the tropical Pl\"ucker
relations! Are the $v(P_i)$ realizable? If so, give a natural linear
space over $\R$, constructed from the $X_i$, that realizes them.
\end{problem}

\section{Series-Parallel Matroids and Linear Spaces} \label{sp}

\begin{figure}
\centerline{\includegraphics{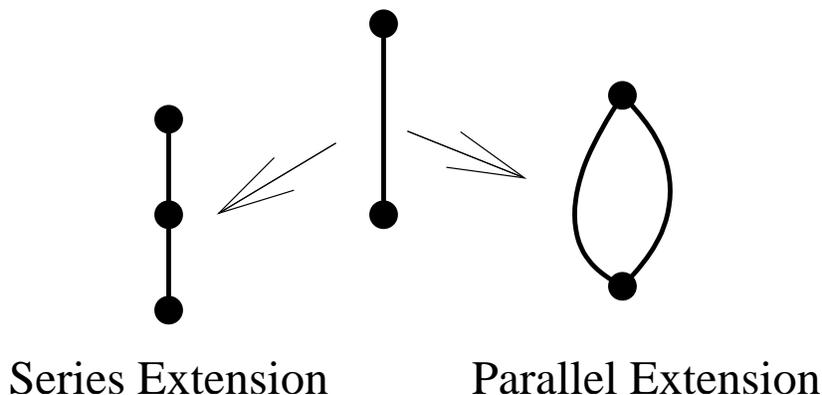}}
\caption{Series and Parallel Extension} \label{SerPar}
\end{figure}

Recall that if $G$ is a connected graph (finite, possibly with
multiple edges or loops) then the associated \emph{graphical matroid}
is a matroid whose elements are the edges of $G$, whose bases are the
spanning trees and whose circuits are the circuits. See \cite{White1},
section 6.1. The rank of this matroid is the number of vertices of $G$
minus one. A matroid is called \emph{series-parallel} if it corresponds
to a series-parallel graph, \emph{i.e.} a graph which can made by
starting from a single edge connecting two distinct vertices by
repeatedly applying the series and parallel extension operators. (See
figure \ref{SerPar}.) For a general introduction to series-parallel
matroids, see Section 6.4 of \cite{White1}.

We will call a matroidal decomposition $\DD$ of $\Delta(s,n)$ a
\emph{series-parallel} decomposition if every facet of $\DD$ is the
polytope of a series-parallel matroid and the associated tropical
linear space will be called a series-parallel tropical linear
space. We will see that series-parallel linear spaces are the most
natural and most manageable tropical linear spaces.

Recall that one of our main aims is to prove.

\begin{f-vector}
Every $d$-dimensional tropical linear space $L$ in $n$-space has at most
$\binom{n-2i}{d-i} \binom{n-i-1}{i-1}$ bounded faces of dimension $i$,
with equality if and only if $L$
is series parallel.
\end{f-vector}

We spend the rest of this section presenting what appears to be a new way
of treating series-parallel matroids, these methods will be applied in
sections \ref{results} and \ref{trees}.

Let $T$ be a trivalent tree (\emph{i.e.} every internal vertex of $T$
has three neighbors) with $n$ leaves whose $n-2$ internal vertices
are colored white and black with $d-1$ black vertices and $n-d-1$
white vertices. For every edge $e$ of $T$ let $A_e \subset [n]$ be the
leaves on one side of $e$ and let $a_e$ be the number of black
vertices on that same side of $e$. Let $\Pi(T)$ be the subpolytope of
$\Delta(d,n)$ consisting of all $(x_i)$ such that, for every edge $e$,
$a_e \leq \sum_{i \in A_e} x_i \leq a_e+1$ .

\begin{prop}
  $\Pi(T)=P_{\mu(T)}$ for a series-parallel matroidal $\mu(T)$. Every series-parallel
matroid can be written as $\mu(T)$ for some bi-colored trivalent tree $T$.

\end{prop}

 When it is necessary to emphasize the dependence of $\Pi(T)$ or $\mu(T)$ on the coloring $c$, we write $\Pi(T,c)$ or $\mu(T,c)$.

\begin{proof}
Our proof is by induction on $n$, when $n=2$ the result is clear. Let $i$
and $j \in [n]$ be two leaves of $T$ that border a common vertex $v$.
Define a new colored trivalent tree $T'$ by deleting $i$ and $j$ from $T$
and forgetting the coloring of $v$, which is now a leaf. By induction,
$\Pi(T')=P_{\mu(T')}$ for a series-parallel matroid $\mu(T')$.

\textbf{Case I:} $v$ is colored white in $T$. Then, cutting $T$ at the
edge separating $i$ and $j$ from the rest of $T$, we see that $x_i+x_j
\leq 1$. We define a map $\phi : \RR^n \to \RR^{n-1}$ by
$\phi(x)_v=x_i+x_j$ and $\phi(x)_k=x_k$ otherwise. An easy inspection of
the defining inequalities shows that $\Pi(T) = \Delta(d,n) \cap\phi^{-1}(\Pi(T'))$. This precisely says that $\Pi(T)$ is the polytope
associated to a series extension of $\mu(T')$ by $v$.

\textbf{Case II:} $v$ is colored black in $T$. This is just like the other
case except that $x_i+x_j \geq 1$, we define $\phi(x)_v=x_i+x_j-1$ and we get a parallel extension.

To show that every series-parallel matroid occurs in this way, reverse the
argument.
\end{proof}

Not all of the inequalities above are necessary to define $\Pi(T)$. The
next
proposition identifies the facets of $\Pi(T)$.

\begin{prop}
The facets of $\Pi(T)$ are given by (1) the equations $\sum_{i \in A_e}
x_e
\geq a_e$ when $e$ is an edge connecting two vertices of the opposite
color and $A_e$ is the set of leaves on the same side of $e$ as the black
endpoint, (2) the equations $x_i \geq 0$ when $i$ is joined to a white
vertex and (3) the equations $x_i \leq 1$ when $i$ is joined to a black
vertex. (All of these are subject to already having the equality $\sum
x_i=d$.)
\end{prop}

\begin{proof}
We first show that all other equations are redundant. Let $e$ be an
internal edge with $A_e$ and $B_e$ the leaves on each side of $e$ and
$a_e$ and $b_e$ the number of black vertices on each side. Then $\sum_{i
\in A_e} x_i =d-\sum_{i \in B_e} x_i$ so $\sum_{i \in A_e} x_i \leq a_e+1$
is equivalent to $\sum_{i \in B_e} x_i \geq d-(a_e+1)=b_e$. If $i$ is a
leaf joined to a white vertex $v$, let $e$ and $e'$ be the other two edges
issuing from $v$ and let $A_e$ and $A_{e'}$ be the sides of $e$ and $e'$
not containing $v$. Then
$$d=x_i+\sum_{j \in A_e} x_j + \sum_{j \in A_{e'}} x_j \geq
x_i+a_e+a_{e'}=x_i+d-1$$
so $x_i \leq 1$. Similarly, if $i$ is a leaf connected to a black vertex,
then the other inequalities imply $x_i \geq 0$.

Suppose that $e$ is an internal edge with endpoints $u$ and $v$, where $u$
is white. Let $e'$ and $e''$ be the other edges containing $u$ and let
$A_e$, $A_{e'}$ and $A_{e''}$ each be as above with $A_e$ on the $u$ side
of $e$ and $A_{e'}$, $A_{e''}$ on the non-$u$ side of $e'$ and $e''$. Then
$$\sum_{i \in A_e} x_e=\sum_{i \in A_{e'}} x_i +\sum_{A_{e''}} x_i \geq
a_{e'}+a_{e''}=a_e.$$
A similar argument show that $\sum_{i \in A_e} x_i \geq a_e$ is redundant
if $v$ is black.

We now must show that each of these inequalities does define a facet. Our
proof is by induction on $n$. First, consider the inequality $\sum_{i \in
A_e} x_e \geq a_e$ where $e$ connects a black vertex $x$ to a white vertex
$y$ and $A_e$ is on the $x$ side of $e$. Let $i$ and $j$ be two leaves of
$T$ that border a common vertex $v$ of $T$. Define $T'$ and $\phi$ as in
the proof of the previous proposition and let $x'$ be a point of $T'$
where all inequalities except for the one arising from the edge $e$ in $T$
are strict. Then any point $x$ obeying $\phi(x)=x'$ and $x_i$, $x_j \in
(0,1)$ will have all inequalities but the required one strict. Similar
arguments apply to the other inequalities.
\end{proof}

We thus see that the facets of $\Pi(T)$ which correspond to loop and
co-loop
free matroids are in bijection with the edges joining vertices of
different colors, the facets corresponding to matroids with a loop
correspond to edges joining leaves to black vertices and the facets
corresponding to matroids with a co-loop correspond to edges from leaves
to white vertices.

Let $T$ be a trivalent tree. Suppose that $e$ is an internal edge. Define
$T
\remove e$ to be the pair of trivalent trees formed as follows: Let $u$
and $v$ be the endpoints of $e$ and let $\{ a, b, v \}$ and $\{ c,d,u \}$
be the neighbors of $u$ and $v$ respectively. Delete the vertices $u$ and
$v$ and all edges ending at them. Then draw new edges joining $a$ to $b$
and $c$ to $d$.  If $e$ is an edge joining a leaf $i$ to a vertex $v$, we
define $T \remove e$ similarly except that we simply delete the vertex
$i$. We will use the same notation $F \remove e$ for the
analogous construction when $F$ is a forest of trivalent trees and $e$ an
internal edge of a tree of $F$. If $T$ is colored, we let $T \remove e$ inherit the coloring of $T$ in the obvious way.

\begin{prop} \label{faces-of-Pi(T)}
Let $e$ be an internal edge of $T$ joining two vertices of opposite
colors. Then the corresponding facet of $\Pi(T)$ is $\Pi(T_1) \times
\Pi(T_2)$ where $T \remove e=T_1 \sqcup T_2$. If $e$ joins a leaf $i$ to
a vertex $v$ then the corresponding facet of $\Pi(T)$ is $\Pi(T \remove e)
\times P_{\{ i \}}$ where $\{ i \}$ is given the structure of a rank $1$
matroid if $v$ is black and a rank $0$ matroid if $v$ is white.

The faces of $\Pi(T)$ are in bijection with the possible colored
forests that can result from repeated splittings of $T$ along edges
connecting vertices of opposite color or connecting leaves to the rest
of $T$. The faces of $\Pi(T)$ that correspond to loop and co-loop-free
matroids are in bijection with the possible colored forests that can arise
by repeated splittings along internal edges between vertices of opposite
color.
\end{prop}

\begin{proof}
The first paragraph follows because it is easy to check that both polytopes are defined
by
the same inequalities. The second paragraph follows by repeatedly using
the
first.
\end{proof}


\section{Special Cases of the $f$-Vector Conjecture} \label{cases}

In this section we will prove the following results:

\begin{thrm}\label{VertexBound}
Let $L$ be a $d$-dimensional tropical linear space in $n$ space. Then
$L$ has at most $\binom{n-2}{d-1}$ vertices, with equality
iff $L$ is series parallel.
\end{thrm}

\begin{thrm}\label{facets}
Let $L$ be a $d$-dimensional tropical linear space in $n$ space with
$n=2d$ or $n=2d+1$. Then $L$ has at most $1$ bounded facet if $n=2d$ and at most $d$ if $n=2d+1$.
\end{thrm}

We first recall the definition of the Tutte polynomial of a
matroid. If $M$ is a rank $d$ matroid on the ground set $[n]$ and $Y \subseteq
[n]$, let $\rho_M(Y)$ denote the rank of $Y$. The polynomial
$$r_M(x,y)=\sum_{Y \subseteq [n]} x^{|Y|-\rho_M(Y)} y^{d -
\rho_M(Y)}$$
is known as the rank generating function of
$M$. The polynomial
$$t_M(z,w)=r_M(z-1,w-1)$$ 
is known as the Tutte polynomial of
$M$. Almost all matroid invariants can be computed in terms of the
Tutte polynomial; see Chapter 6 of \cite{White2} for a survey of its
importance. Write $t_M(z,w) = \sum t_{ij} z^i w^j$. Although not
obvious from this definition, all of the $t_{ij}$ are nonnegative. For
$n \geq 2$ we have $t_{10}=t_{01}$, this number is known as the beta
invariant of $M$ and denoted $\beta(M)$.

We will need the following result:

\begin{prop}\label{Crapo}
Let $M$ be a matroid on at least $2$ elements. Then $\beta(M)=0$ if and
only if $M$ is disconnected and $\beta(M)=1$ if and only if $M$ is
series-parallel.
\end{prop}

\begin{proof}
See \cite{Crapo}, theorem II, for the first statement and \cite{Bryl},
theorem 7.6, for the second.
\end{proof}

The key to proving theorem \ref{VertexBound} will be proving the following
formula:

\begin{Lemma}\label{TutteLemma}
Let $M$ be a matroid and let $\DD$ be a matroidal subdivision of
$P_M$. Let $\oDD$ denote the set of internal faces of $\DD$. Then
$$t_M(z,w) = \sum_{P_\gamma \in \oDD} (-1)^{\dim(P_M) - \dim(P_\gamma)}
t_{\gamma}(z,w).$$
\end{Lemma}

Before proving lemma \ref{TutteLemma}, let us see why it implies
theorem \ref{VertexBound}. Considering the case where $M=\Delta(d,n)$
and comparing the coefficients of $z$ on each side, we see that
$$\beta(\Delta(d,n))=\binom{n-2}{d-1}=\sum_{\substack{ \gamma \in
\DD \\
\gamma \textrm{\ a\ facet} }} \beta(\gamma).$$
(Note that all nonfacets of $\oDD$ correspond to disconnected matroids and
hence have beta invariant $0$. and that every facet of $\DD$ is in
$\oDD$.)

Every term in the right hand sum is a positive integer. Thus, $\DD$ has at
most  $\binom{n-2}{d-1}$ facets, with equality if and only if all of the
facets have beta invariant $1$, \emph{i.e.}, if and only if $\DD$ is a
decomposition into series-parallel matroids. \qedsymbol

We now prove lemma \ref{TutteLemma}.

\begin{proof}

Since $t_M(z,w)=r_M(z-1,w-1)$, it is enough to prove
$$r_M(x,y) =
\sum_{P_\gamma \in \oDD} (-1)^{\dim(P_M) - \dim(P_\gamma)}
r_{\gamma}(x,y).$$
Plugging in the definition of $r_M$ and interchanging summation, it is
enough to show that for every $Y \subseteq [n]$,
$$x^{|Y|-\rho_M(Y)} y^{d - \rho_M(Y)}=\sum_{\gamma \in \oDD}
(-1)^{\dim(P_M) - \dim(P_\gamma)} x^{|Y|-\rho_\gamma(Y)}
y^{d -
\rho_\gamma(Y)}.$$

Comparing coefficients of
$x^{|Y|-r} y^{d-r}$, we are thus being asked to show that
$$\sum_{\substack{P_\gamma \in \oDD \\ \rho_{\gamma}(Y)=r}} (-1)^{\dim(P_M)
-
\dim(P_\gamma)} = \left\{ \begin{matrix} 1 \hbox{\ if\ $r=\rho_M(Y)$} \\ 0
\hbox{\ if\ $r<\rho_M(Y)$} \end{matrix} \right\} . $$
The sum is empty if $r>\rho_M(Y)$. Equivalently, we will show
$$\sum_{\substack{P_\gamma \in \oDD \\ \rho_{\gamma}(Y) \geq r}}
(-1)^{\dim(P_M) - \dim(P_\gamma)} = 1$$
for all $r \leq \rho_M(Y)$.

Let $\ell_Y$ be the linear function $\Delta(d, n) \to \RR$ mapping
$(x_i) \mapsto \sum_{i \in Y} x_i$. Then $\rho_{\gamma}(Y)=\max_{x \in
P_\gamma} \ell_Y(x)$. Thus, we see that $\rho_{\gamma}(Y) \geq r$ if
and only if $\gamma$ has a nonempty intersection with the half space
$\ell_Y > r-1/2$.  The promised equality now follows by the following
lemma applied to the polytope $P_M \cap \{ x : \ell_Y(x) > r-1/2 \}$.

\begin{Lemma}
Let $P$ be any bounded polytope and $\Gamma$ the internal faces of a
decomposition of $P$. Then $\sum_{\gamma \in \Gamma} (-1)^{\dim(P) -
\dim(\gamma)}=1$.
\end{Lemma}

\begin{proof}
This sum is $(-1)^{\dim P} (\chi(P) - \chi(\partial P))$ where $\chi$ is
the Euler characteristic. As $P$ is contractible and $\partial P$ is
a sphere of dimension ${\dim(P)-1}$, the result follows.
\end{proof}

\end{proof}

\begin{figure}
\centerline{\includegraphics{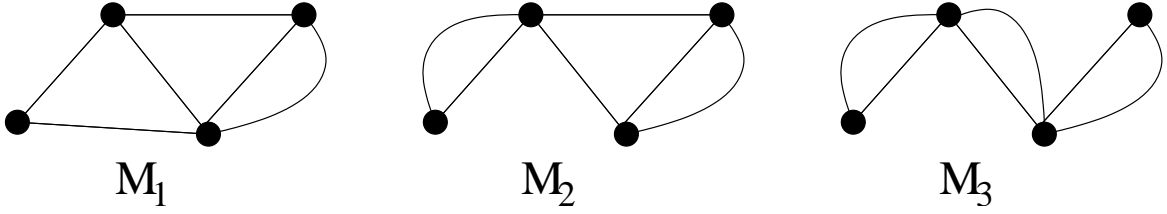}}
\caption{The Graphical Matroids $M_1$, $M_2$ and $M_3$} \label{BadTutte}
\end{figure}

One might hope to use the higher degree terms of lemma \ref{TutteLemma} to
produce additional bounds on the $f$-vector of $\DD$. Unfortunately, lemma
\ref{TutteLemma} is incapable of producing a complete set of restrictions.
For example, consider the matroids $M_1$, $M_2$ and $M_3$ corresponding to
the graphs in figure \ref{BadTutte}.

The Tutte polynomials of these matroids are

\begin{alignat*}{4}
t_{M_1} =& z+w+ & 2z^2+3zw+2w^2+ & z^3+\phantom{2}z^2 w+\phantom{2}w z^2+w^3 \\
t_{M_2} =&      &  \phantom{2}z^2+2zw+\phantom{2}w^2 + & z^3+2z^2 w+2z w^2+w^3 \\
t_{M_3} =&      &                & z^3+3z^2 w+3z w^2+w^3.
\end{alignat*}

We have
$$6t_{M_1}-9t_{M_2}+4t_{M_3}=t_{\Delta(3,6)}=6z+6w+3z^2+3w^2+z^3+w^3.$$
Nonetheless, it follows from theorem \ref{facets} (proved below) that
there is no matroidal decomposition of $\Delta(3,6)$
which has $4$ internal faces of codimension $2$.

\begin{proof}[Proof of theorem \ref{facets}]

Suppose that $L$ is a $d$-dimensional tropical linear space in $n$-space.
Let $M^{\vee}$ be a bounded facet of $L$, dual to a polytope $P_M$. Then
the matroid $M$ has rank $d$ and has $d$ connected components. So
$M=\bigoplus_{i=1}^d M_i$ where each $M_i$ has rank $1$ and, because
$M^{\vee}$ is bounded, $|M_i| \geq 2$. Let $S_i$ be the ground set of
$M_i$, the vertices of $P_M$ are those subsets of $[n]$ that meet each
$S_i$ precisely once.

Now, consider two such facets $M^{\vee}$ and $(M')^{\vee}$ occuring in the
same linear space and define $S_i$ and $S'_i$ as above. The polytope $P_M
\cap P_{M'}$ must be either matroidal or empty.

\begin{lemma}
Let $\Gamma=\Gamma(M,M')$. $\Gamma$ may have multiple edges, let $\overline{\Gamma}$
be the simple graph obtained by replacing each multiple edge of $\Gamma$
by a single edge. Then $P_M \cap P_{M'}$ is empty if and only if
$\overline{\Gamma}$ has no perfect matchings. $P_M \cap P_{M'}$ is
matroidal if and only if $\overline{\Gamma}$ has precisely one perfect
matching.
\end{lemma}

A perfect matching of a graph $\Gamma$ is a collection of edges such that, for every vertex $v$ of $\Gamma$, exactly one edge in the collection contains $v$.

\begin{proof}
The vertices of $P_M \cap P_{M'}$ are in bijection with the perfect
matchings of $\Gamma$. There is an obvious surjection from perfect
matchings of $\Gamma$ to matchings of $\overline{\Gamma}$. This proves the
first claim.

Suppose that $P_M \cap P_{M'}$ is nonempty and matroidal. Let
$(B_1,B_2)$ be an edge of $P_M \cap P_{M'}$, then $B_1$ and $B_2$
correspond to perfect matchings of $\Gamma$ that differ only by a
single edge say $B_1 \setminus \{ e_1 \} = B_2 \setminus \{ e_2
\}$. Then $e_1$ and $e_2$ have the same image in $\overline{\Gamma}$
in order for both $B_1$ and $B_2$ to be perfect matchings. Since every
pair of vertices of $P_M \cap P_{M'}$ is connected by a chain of
edges, we see that every perfect matching of $\Gamma$ gives rise to
the same perfect matching of $\overline{\Gamma}$ or, in other words,
$\overline{\Gamma}$ has exactly one perfect matching. Conversely, if
$\overline{\Gamma}$ has exactly one perfect matching, it is clear that
$P_M \cap P_{M'}$ is a product of simplices and hence matroidal.
\end{proof}

We now turn to the cases we are interested in, when $n=2d$ or $2d+1$. We maintain the
notation of the lemma. First, suppose that $n=2d$. Then every $S_i$ must
have order $2$ so, for any $S_i$ and $S'_i$, $\Gamma$ must consist of
disjoint cycles. The only way that $\overline{\Gamma}$ can be a forest is
if all of those cycles have length $2$. But then $\{S_i \}=\{ T_i \}$. So
we see that a linear space of dimension $d$ in $2n$ space can have at most
one bounded facet.

Now suppose that $n=2d+1$.

Then all of the $S_i$ have order $2$ except
for one which has order $3$. We adopt the convention that $|S_1|=3$.
$\Gamma$ consists of several disjoint cycles and one component $C$. $C$
consists of a pair of vertices $v$ and $w$ of degree $3$ with opposite
colors and three disjoint paths $\gamma_1$, $\gamma_2$, $\gamma_3$.
Either these paths all run from $v$ to $w$ or there is one path each running from
$v$ to $v$, from $v$ to $w$ and from $w$ to $w$.

In any of these cases, $\Gamma$ has at least one perfect matching. One
can check that the only case where $\overline{\Gamma}$ has only one
perfect matching is when all of the cycles of $\Gamma \setminus C$ have
length $2$ and $C$ consists of a path from $v$ to $w$ and two cycles of
length $2$, one containing $v$ and the other containing $w$. The second condition is equivalent to requiring that $S'_{i'} \subset S_1$ and $S_i \subset S'_1$ for some $i$ and $i'$.

Now, consider all of the three element subsets of $[n]$ that occur as an
$S_1$ for some bounded facet of $L$. We first note that, if $S_1=T_1$ then
$S_i=T_i$ for every $i$ (after reordering) as otherwise $\Gamma$ would
have a cycle of length more than $2$. So all of these three element
subsets are different. Also, if $|S_1 \cap T_1|=2$ then there are two
paths between $S_1$ and $T_1$ and $\overline{\Gamma}$ has more than one
perfect matching.

Let $f$ be the number of bounded facets of $L$. Build a graph whose
vertices are the various $S^r_1$ and the elements $\{ 1, \ldots, n \}$ and
where there is an edge between $S^r_1$ and $k$ if $k \in S^r_1$. This graph
has $f+n=f+2d+1$ vertices and $3f$ edges. Thus, if $f>d$, this graph has
as many edges as vertices and must contain a cycle. In other words, there
exist $i_1$, \dots, $i_f$, $j_1$, \dots, $j_f \in [n]$ and facets
corresponding to collections $(S^r_s)$, $1 \leq s \leq d$, $1 \leq r \leq
f$ such that $S^r_1 = \{ i_r, i_{r+1}, j_r \}$, where the index $r$ is
cyclic modulo $f$.

Now, for any $r$ and $r'$, there must be some $S^{r}_s \subset S^{r'}_1$. Consider first the case where
$r'=r+1$. Then $\{ j_{r+1}, i_{r+2} \}$ must be one of the $S^{r}_s$, as
we already know $i_{r+1} \in S^{r}_1$. Now consider when $r'=r+2$. We get
that $\{ j_{r+2}, i_{r+3} \}$ must be one of the $S^r_{s}$, as $i_{r+2}$,
the other member of $S^{r'}_1$, is already contained in  $\{ j_{r+1},
i_{r+2} \}$. Continuing in this manner, we get that $\{ j_{r+k}, i_{r+k+1}
\}$ is an $S^r_s$ for every $k$. But, when $k=f-1$, this intersects
$S^r_1$, contradicting that the $S^r_i$ are distinct.

\end{proof}

\begin{remark}
The $f$-vector conjecture predicts that, if $n=2d+e$,
$L$ should have at most $\binom{d+e-1}{e}$ bounded facets. In the
cases $e=0$ and $e=1$, $\Gamma$ always had at least one matching and the
requirement that $\overline{\Gamma}$ have exactly one is very
restrictive. Once $e=2$, it is possible for $\overline{\Gamma}$ to have no perfect matchings and the problem grows much harder.

Nonetheless, I suspect that this kind of reasoning, looking only at
the facets of $L$ and not at the smaller spaces that contain them, is
adequate to resolve the $f$-vector conjecture for facets once $e
\geq 2$. More specifically, I conjecture:
\end{remark}

\begin{problem}
Let $n=2d+e$. It is not possible to find more than $\binom{d+e-1}{e}$ distinct
partitions of $[n]$ into $d$ disjoint subsets, $[n]=S^r_1 \sqcup \cdots
\sqcup S^r_{d}$ where each $S^r_s$ has order at least $2$ and, for every
$r$ and $r'$, the graph $\overline{\Gamma}$ formed from $\{ S^r_1, \ldots,
S^r_d \}$ and $ \{ S^{r'}_1, \cdots, S^{r'}_d \}$ has at most one perfect
matching.
\end{problem}

\section{Results on Constructible Spaces}
\label{results}

In this section, we will prove several of the previously stated results on contructible spaces.

Our first aim is to prove theorem \ref{cons-sp}: that every
constructible space is series-parallel. Clearly, $L^{\perp}$ is
series-parallel if $L$ is, since the dual of a series-parallel matroid is
series-parallel. So it is enough to show that, if $L$ and $L'$ are two
series-parallel tropical linear spaces that meet transversely then $L
\cap L'$ is as well.

Let $L$ and $L'$ be two tropical linear spaces of dimensions $d$ and $d'$
in $n$ space which meet transversely at a point $w$. Let $\DD$ and $\DD'$
be the appropriate subdivisions of
$\Delta(d,n)$ and $\Delta(d',n)$. Let $M^{\vee}$
and $(M')^{\vee}$ be the faces of $\DD^{\vee}$ and $(\DD')^{\vee}$
containing $w$. Let $M=\bigoplus_{i \in I} M_i$ and $M'=\bigoplus_{i \in
I'} M'_{i}$ be the connected components of $M$ and $M'$. We claim that all of the $M_i$ and $M'_i$ are series parallel. $P_M$ will be a face of some facet $P_{\tilde{M}}$ of $\DD$, so this follows from

\begin{Lemma}
Let $\tilde{M}$ be a series-parallel matroid and let $P_M$ be a face
of $P_{\tilde{M}}$, assume that $M$ is loop-free. Let $M=\bigoplus
M_i$ be the decomposition of $M$ into connected components. Then each
of the $M_i$ are series-parallel.
\end{Lemma}

\begin{proof}
This follows from the description of the faces of $\tilde{M}$ in lemma
\ref{faces-of-Pi(T)}. (Note that a rank one matroid on a single
element is series-parallel.)
\end{proof}

The result now follows from
\begin{prop}
Let $M=\bigoplus M_i$ and $M'=\bigoplus M'_i$ be transverse matroids on $[n]$ with all of the $M_i$ and $M'_i$ series-parallel. Then $M \wedge M'$ is a direct sum of series-parallel matroids, and is series parallel if it is connected.
\end{prop}

\begin{proof}
$M \wedge M'$ is formed by a sequence of parallel connections. The parallel connection of two series-parallel matroids is series-parallel.
\end{proof}

Our next goal is to prove theorem \ref{cons-f} -- every constructible space achieves the $f$-vector of the $f$-vector
conjecture.

\begin{Lemma} \label{symmetry}
Let $M$ be a series-parallel matroid. Call a flat $Q$ of $M$ ``good'' if
$M|_Q$ and $M/Q$ are connected and not loops. Let $f$ and $g \in M$. Then
the collection
of all good flats $Q$ such that $f \in Q$, $g \not\in Q$ forms a chain
$Q_1 \subset \cdots \subset Q_d$. Moreover, the length $d$ of this chain
is
preserved when the roles of $f$ and $g$ are switched.
\end{Lemma}

We define $d(f,g)$ to be $d$.

\begin{proof}
If $Q$ is good then $M|_Q \oplus M/Q$ corresponds to a facet of $M$. Let
$M=\mu(T)$. By lemma \ref{faces-of-Pi(T)}, all of the facets of $P_M$
correspond to either edges of $T$ which
contain a leaf or edges of $T$ that connect to vertices of opposite
colors. Under that correspondence, an edge $e$ corresponds to a good flat
that contains $f$ and not $g$ if and only if $e$ separates $f$ and $g$
with the end closer to $g$ colored black and the end closer to $f$
either colored white or equal to $e$. Consider the path $\gamma$ through
$T$ connecting $f$ and $g$, divide $\gamma$ into alternating blocks of
white and black vertices. Then we see that the number of good flats
containing $f$ and not $g$ is equal to the number of black blocks. This
number is the same when $f$ and $g$ are switched. The flat corresponding
to $e$ consists of the leaves of $T$ that are on the non-black end of $e$.
Thus, all of the flats in question are nested.
\end{proof}

Our first major goal on the way to proving theorem \ref{cons-f}, nontrivial in its own right, is to prove:

\begin{thrm} \label{translate}
Let $L$ and $L'$ be two series-parallel tropical linear spaces in
$n$-space, of dimensions $d$ and $d'$. Let $a$ and $b \in \RR^n$ be such
that $L$ meets $L'+a$ and $L'+b$ transversely. Then $L \cap (L'+a)$ and $L
\cap (L' + b)$ have the same bounded $f$-vector.
\end{thrm}

\begin{proof}[Proof of theorem \ref{translate}]
Let $U_i \in \RR^n$ be the set of $v$ such that, for any two faces $F$
and $F'$ of $L$ and $L'$, either the realtive interiors of $F$ and $F'+v$
are
disjoint or they span an affine linear space of dimension at least
$n-i$. Thus, $L$ and $L'+v$ meet transversely if and only if $v \in
U_0$. $\RR^n \setminus U_i$ is a polyhedral complex of dimension at
most $n-i-1$. In particular, we see that $U_1$ is path connected. We may
therefore join $a$ and $b$ by a path through $U_1$. It is thus enough
to show that, if $v \in U_1 \setminus U_0$ and we perturb $v$ into
$U_0$ then the $f$-vector we get is independent of the choice of
perturbation. Without loss of generality, we may assume that $v=0$.

Let $N^{\vee}$ be a face of $L \scap L^{\vee}$ at which two faces
$M^{\vee}$ and $(M')^{\vee}$ fail to meet transversely. Let the
connected components of $M$ and $M'$ be $\bigoplus_{i \in I} M_i$ and
$\bigoplus_{i \in I'} M'_i$ and let $\Gamma=\Gamma(M,M')$.

This time, since
$M^{\vee}$ and $(M')^{\vee}$ only span an affine hyperplane, $\Gamma$
will have first betti number $1$ and will thus contain a unique
cycle. Let $e_1$, \dots, $e_{2r}$ be the edges of that cycle in cyclic
order. ($\Gamma$ is bipartite so the cycle has even length.) 
Write $v_i$ and $v_{i+1}$ for the ends of $e_i$ and write
$M^{(')}_{v}$ to mean $M_v$ if $v \in I$ and $M'_v$ if $v \in I'$. 

Now, let $w \in \RR^n$ and let $\EE_w$ be the decomposition of
$\Delta(d+d'-n,n)$ associated to $L \scap (L'+w)$. For $\epsilon$ a
sufficiently small positive real number, $\EE_{ew}$ is a subdivision of
$\EE_0$. We will explore how the interior of $P_N$ decomposes in
$\EE_{\epsilon w}$.

\begin{Lemma}\label{perturb}
There exist decompositions $\EE_+$ and $\EE_-$ of $P_N$ such that, if $\pm
\sum_{i=1}^{2r} (-1)^i w_{e_i}>0$ then $N$ is subdivided according to
$\EE_{\pm}$ in $\EE_{\epsilon w}$ for $\epsilon$ a sufficiently small
positive number.  Moreover, there is an integer $D$ such that $\EE_+$ and
$\EE_{-}$ each have $\binom{D}{c+1}$ interior faces of codimension $c$.
\end{Lemma}

The theorem follows from this lemma, as any generic perturbation of $L
\scap L'$ to $L \cap (L'+w)$ is going to have $\sum (-1)^i w_{e_i} \neq
0$ and this lemma shows that the $f$-vector changes in the same way under
all such perturbations. (There may be several different $P_N$ which all
subdivide under perturbation, but we may apply the lemma to each of them.)

\end{proof}

\begin{proof}[Proof of lemma \ref{perturb}:]
By assumption, $M^{\vee}$ and $(M')^{\vee}$ together span a
hyperplane, call this hyperplane $H$. Explicitly, the equation of $H$
is $\sum (-1)^i x_{e_i}=0$. If we perturb $L'$ to $L'+w$ for $w
\not\in H$, $M^{\vee}$ and $(M')^{\vee}+w$ will not meet. Instead,
there will be faces $P_F$ and $P_{F'}$ of $P_M$ and $P_{M'}$ with
$M^{\vee} \subset F^{\vee} \subset L$ and $(M')^{\vee} \subset
(F')^{\vee} \subset L'$, such that the relative interiors of
$F^{\vee}$ and $(F')^{\vee}$ were previously disjoint but now meet and
contribute faces to the subdivision of $P_N$. For
any given pair $(F,F')$, whether the relative interiors of $F^{\vee}$
and $(F')^{\vee}+\epsilon w$ meet or not for small $\epsilon$ will
depend only on the side of $H$ on which $w$ lies.

The facets occurring in the subdivision of $P_N$ will came from pairs
$(F,F')$ where either $F=M$ and $F'$ is a loop-free facet of $M'$ or
\emph{vice versa}. We will concentrate on the former case. Every
loop-free facet of $M$ is of the form $M|_{Q} \oplus M/Q$ for $Q$ a
good flat of some particular irreducible component $M_i$ of $M$. The
local cone of $F^{\vee}$ at $0$ is the span of the linear space
spanned by $M^{\vee}$ and a ray in the direction $\sum_{q \in Q}
e_q$. Note that $Q \cap \{ e_1, \ldots, e_{2r} \}$ contains at most
$2$ elements. We thus see that the image of this cone in the one
dimensional space $\RR^n/H$ is $0$ unless exactly one of $\{ e_1,
\ldots, e_{2r} \}$ lies in $Q$. If this cone maps to $0$ in
$\RR^n/H$ then $F^{\vee} \cap ((F')^{\vee}+w)=\emptyset$ for $w
\not\in H$. If $Q \cap \{ e_1, \ldots, e_{2r} \}=\{ e_k \}$ then $F$
lies on one side of $H$ if $k$ is even and the other if $k$ is
odd. So, if $\sum_{i=1}^{2r} (-1)^i e_{w_i}>0$, we have $F^{\vee} \cap
((F')^{\vee}+w)=\emptyset$ if $k$ is even and is nonempty if $k$ is
odd; the reverse holds if $\sum_{i=1}^{2r} (-1)^i e_{w_i}<0$.

After running through the same argument in the case that $F'=M'$ and $F$
is a facet of $M$, we see that, if
$\sum (-1)^i w_{e_i}>0$, then $P_N$ is divided into a number of facets
equal to the number of good flats containing $e_i$ and not $e_{i+1}$ as
$i$ runs from $1$ to $2d$. In other words, $P_N$ is divided into
$\sum_{i=1}^{2r} d(e_i, e_{i+1})$ facets. When $\sum (-1)^i w_{e_i}<0$,
$P_N$ is divided into $\sum_{i=1}^{2r}
d(e_{i+1},e_i)$ facets. These sums are equal to a common value, call it
$D$, by lemma \ref{symmetry}.

 A
similar argument shows that the number of codimension $c$ faces when $\sum
(-1)^i w_{e_i}>0$  is the
number of ways of, for each $1 \leq i \leq 2r$, choosing a chain
$Q^i_1 \subset \cdots \subset \cdots Q^i_{q_i}$ of flats in
$M^{(')}_{v_i}$
containing $e_{i-1}$ and not $e_i$ subject to the condition $\sum
q_i=c+1$. This is clearly the same as choosing $c+1$ of the $D$ flats
counted above, so there are $\binom{D}{c+1}$ codimension $c$ faces. A
practically identical argument shows the same in the case  $\sum (-1)^i
w_{e_i}<0$.
\end{proof}

We now are ready to begin our final assault on theorem
\ref{cons-f}. We will first need to study how the various faces of
$\Delta(d,n)$ are subdivided in a constructible decomposition. Let $S$
and $T$ be disjoint subsets of $[n]$. Write $\Delta(d,n) \setminus 
S/T$ for the face of $\Delta(d,n)$ given by the equations $x_i=0$ for
$i \in S$ and $x_i=1$ for $i \in T$. If
$p$ is a tropical Pl\"ucker vector then define $p \setminus  S/T$ to
be the function on the vertices of $\Delta(d-|T|, n-|S|-|T|)$ given by
$(p \setminus  S/T)_I=p_{I \cup T}$. It is clear that $p \setminus
 S/T$ obeys the tropical Pl\"ucker relations, as they are a subset
of the tropical Pl\"ucker relations for $p$. So $p \setminus  S /T$
corresponds to a tropical $d-|T|$ plane in $n-|S|-|T|$ space. If $L$
is the tropical linear space corresponding to $p$ then we write $L
\setminus  S/T$ for the tropical linear space corresponding to $p
\setminus  S /T$.

\begin{lemma}
If $L$ is a constructible $d$-plane in $n$-space and $S$ and $T$ are
disjoint subsets of $[n]$ then $L\setminus  S/T$ is also constructible.
\end{lemma}

\begin{proof}
If $L$ is a hyperplane then every $L \setminus  S/T$ is either a
hyperplane
or the whole space in which it sits, and hence is constructible.  Now suppose
that the lemma holds for $L$. Then it holds for $L^\perp$ as $L^\perp
\setminus 
S/T=(L \setminus  T/S)^{\perp}$. Suppose that the lemma holds for $L$
and
$L'$.
Then it holds for $L \scap L'$ as $(L \scap L') \setminus  S/T=(L
\setminus 
\emptyset / (S \cup T)) \scap (L' \setminus  S/T)$.
\end{proof}

\begin{proof}[Proof of theorem~\ref{cons-f}]
Our proof is by induction, first on $n$ and then on the number of steps
used to construct $L$. We set
$$f_{i,d,n}=\binom{n-2i}{d-i} \binom{n-i-1}{i-1}=\frac{(n-i-1)!}{(i-1)!
(d-i)! (n-d-i)!}.$$

First, suppose that the theorem is true for $L$. $L^{\perp}$ and $L$ have
the
same bounded $f$-vector, so the theorem for $L^{\perp}$ follows from the
symmetry
$f_{i,d,n}=f_{i,n-d,n}$. For the rest of the proof, suppose that $L$ and
$L''$ are a constructible $d$-plane and $d'$ plane in $n$ space which meet
transversely and for which we know the theorem to hold. In addition, we
know that the theorem holds for every $L \setminus  S/T$ and $L''
\setminus 
S/T$ by our inductive hypothesis. Our aim will be to prove the theorem for
$L \cap L''$.

By theorem \ref{translate}, it is enough to prove the result for $L \cap
(L''+w)$ where $w$ is chosen generically enough. We will take a $w$ for
which $w_1 \ll  w_2 \ll  \cdots \ll  w_n$ and otherwise generic. Set
$L'=L''+w$.

Suppose that $M^{\vee}$ is a face of $L$ dual to $P_M \subset
\Delta(d,n)$. Then $M$ is loop-free, let $T$ be the set of isthmi of
$M$ so that $M$ lies in the relative interior of $\Delta(d,n) \setminus 
\emptyset/T$. Then $M^{\vee}/(1,\ldots, 1)$ is the product of
$\RR_{+}^{T}$ and a bounded polytope. Let $(M'')^{\vee}$ be a face of
$L''$ and define $T'$ similarly. We write $(M')^{\vee}=(M'')^{\vee}+w$ and
$M'=M''$.

Define an ordered pair $(T, T')$ of subsets of $[n]$ to be \emph{nice} if,
for all $1 \leq i<j \leq n$, either $j \in T$ or $i \in T'$.

\begin{lemma}
For $w_1 \ll  w_2 \ll  \cdots \ll  w_n$, we have $M^{\vee} \cap
(M')^{\vee} \neq
\emptyset$ if and only if $(T,T')$ is nice. If $(T,T')$ is nice then
$M^{\vee} \cap (M')^{\vee}$ is bounded if and only if $T \cap
T'=\emptyset$.
\end{lemma}.

\begin{proof}
First, suppose that $(T,T')$ is not nice, so there exists some $i<j$ with $j
\not \in T$ and $i \not \in T'$. On $M^{\vee}$, $x_j-x_i$ is bounded above
(as $j \not \in T$). On $(M'')^{\vee}$, $x_j-x_i$ is similarly bounded
below. The value of $x_j-x_i$ on $(M'')^{\vee}$ is $w_j-w_i$ larger than
on $(M'')^{\vee}$, so if $w_j-w_i$ is large enough then the values of this
functional on  $M^{\vee}$ and $(M'')^{\vee}+w$ will be distinct.  Thus,
$M^{\vee}$ and $(M'')^{\vee}+w$ will be disjoint.

Now suppose that $(T,T')$ is nice. Let $(x_i) \in M^{\vee}$ and $(x''_i)
\in (M'')^{\vee}$, set $x'_i=x''_i+w_i$ and set $d_i=x'_i-x_i$. For $w$
large enough, we have $d_1 < d_2 < \cdots < d_n$. The hypothesis that
$(T,T')$ is good implies that there is some $b \in [n]$ such that $i \in
T'$ for every $i <b$ and $j \in T'$ for every $j >b$. After translating
$(x_i)$ by $(1,\ldots,1)$, which will not change the assumption that
$(x_i) \in M$, we may assume that $x_b=x'_b$. Consider the point
$$z=(x_1,\ldots, x_{b-1}, x_b, x_{b+1}+d_{b+1}, \ldots x_n+d_n) =
(x'_1-d_1,\ldots, x'_{b-1}-d_{b-1}, x'_b, x'_{b+1}, \ldots x'_n).$$
Clearly, $z \in (x_i)+\RR_{\geq 0}^T$ and  $z \in (x'_i)+\RR_{\geq
0}^{T'}$, so $z \in M^{\vee} \cap (M')^{\vee}$. We have now proved the
first claim.

Now assume that $(T,T')$ is nice, so
there is some point $z=(z_i) \in M^{\vee} \cap (M')^{\vee}$. Suppose that
$i \in T \cap T'$. Then $z+u e_i \in M^{\vee} \cap (M')^{\vee}$ for every
$u>0$ and $M^{\vee} \cap (M')^{\vee}$ is unbounded.

Suppose, on the other hand, that $T \cap T'=\emptyset$. Consider any
$i<j$. If neither $i$ nor $j$ is in $T$ then $x_i-x_j$ is bounded on
$M^{\vee}$; if neither $i$ nor $j$ is in $T'$ then $x_i-x_j$ is bounded
on $(M')^{\vee}$; if $j \in T$ and $i \in T'$ then
$x_i-x_j$ is bounded above on $M^{\vee}$ and  below on $(M')^{\vee}$. In
any of these cases, we see that $x_i-x_j$ is bounded on $M^{\vee} \cap
(M')^{\vee}$. So all of the $x_i-x_j$ are bounded on $M^{\vee} \cap
(M')^{\vee}$ which implies that $(M^{\vee} \cap (M')^{\vee})/(1,\ldots,1)$
is bounded.
\end{proof}

$M^{\vee}$ corresponds to a bounded face of $L \setminus  \emptyset/T$.
Writing $t=|T|$, if this bounded face is $j$-dimensional then $M^{\vee}$
is $(j+t)$-dimensional. Similarly writing $t'=|T'|$ and $j'$ for the
dimension of the face of $L' \setminus  \emptyset /T'$ corresponding to
$(M')^{\vee}$, the dimension of $(M')^{\vee}$ is $j'+t'$. The dimension of
$M^{\vee} \cap (M')^{\vee}$ is then $j+j'+t+t'-n$. We thus see that the
number of $i$-dimensional bounded faces of $L \cap L'$ is
$$\sum_{\substack{ (T,T') \textrm{\ nice } \\ T \cap T'=\emptyset }}
\sum_{j+j'+t+t'-n=i} f_{j,d-t,n-t} f_{j',d'-t',n-t'}.$$

The only nice pairs $(T,T')$ with $T \cap T'=\emptyset$ are $( \{ 1,
\ldots, b-1 \}, \{ b+1, \ldots, n \})$ and  $( \{ 1, \ldots, b \}, \{ b+1,
\ldots, n \})$. So we can write this sum as
$$\sum_b \left(\sum_{j+j'=i+1} f_{j, d-b+1, n-b+1} f_{j', d'-n+b, b} +
\sum_{j+j'=i} f_{j,d-b,n-b} f_{j', d'-n+b, b} \right).$$
We will now show that this sum is $f_{i,D,n}$, where $D$ is $d+d'-n$.

We work with the first sum first and we use the formula
$f_{i,d,n}=\binom{n-d-1}{i-1} \binom{n-i-1}{n-d-1}$. Exchanging order of
summation, we have
$$\sum_{j+j'=i+1} \binom{n-d-1}{j-1}   \binom{n-d'-1}{j'-1}     \sum_b
\binom{n-b-j}{n-d-1} \binom{b-j'-1}{n-d'-1}.$$

Using the identity $\sum_{k+k'=p} \binom{k}{q}
\binom{k'}{r}=\binom{p+1}{q+r+1}$, the inner sum evaluates to
$\binom{n-i-1}{2n-d-d'-1}=\binom{n-i-1}{n-D-1}$. We are reduced to
$$\binom{n-i-1}{n-D-1} \sum_{j+j'=i+1} \binom{n-d-1}{j-1}
\binom{n-d'-1}{j'-1}.$$
The remaining sum can be evaluated by the identity $\sum_{j+j'=r}
\binom{p}{j} \binom{q}{j'}=\binom{p+q}{r}$. We get $\binom{n-i-1}{D-i}
\binom{n-D-2}{i-1}$.

A similar argument shows that the second sum is $\binom{n-i-1}{D-i}
\binom{n-D-2}{i-2}$. Adding these together, we get
$$\binom{n-i-1}{D-i} \left[ \binom{n-D-2}{i-1} +  \binom{n-D-2}{i-2}
\right]
=
\binom{n-i-1}{D-i}  \binom{n-D-1}{i-1}=f_{i,D,n}.$$
\end{proof}

\section{Tree Tropical Linear Spaces} \label{trees}

In this section, we will describe some tropical linear spaces which
achieve the bounds of the $f$-vector conjecture and have very elegant
combinatorics. First, we will need to summarize some results from
\cite{SpeySturm} concerning the case $d=2$.

Let $T$ be a tree with leaves labeled by $[n]$ and no vertices of
degree $2$. Let each edge $e$ of $T$ be assigned a length $\ell(e)$,
which is a positive real number if $e$ is internal (\emph{i.e.} if
neither endpoints of $e$ is a leaf) and can be any real number if one
of $e$'s endpoints is a leaf. For $i$ and $j \in [n]$, set
$d(i,j)=\sum \ell(e)$ where the sum is over all edges $e$ separating
$i$ and $j$.

\begin{thrm}
Let $L$ be a tropical $2$-plane in $n$-space. Then there exists a tree $T$
as above such that $p_{ij}=-(1/2)d(i,j)$. Moreover, $L/(1,\ldots,1)$ is an
embedding of $T \setminus \{ \hbox{\emph{leaves\ of\ }}T \}$ in
$\RR^n/(1,\ldots,1)$, where the bounded part of $L/(1,\ldots,1)$ is the
image of the internal edges of $T$. For every such $T$, set $p_{ij}=(-1/2)
d(i,j)$, $p_{ij}$ obeys the tropical Pl\"ucker relations. $p$ is generic
if and only if the corresponding $L$ is series-parallel if and only if the
tree $T$ is trivalent (\emph{i.e.} every non-leaf has degree $3$.) Every
$2$-dimensional tropical linear space is realizable.
\end{thrm}

\begin{proof}
See \cite{SpeySturm}, section 4.
\end{proof}

Note that this implies that $L/(1,\ldots1)$ has at most $n-2$ vertices and
at most $n-3$ bounded edges, with equality precisely when $L$ is
series-parallel. So the $f$-vector conjecture is true when $d=2$.

Now, suppose that $p_{ij}$ obeys the tropical Pl\"ucker relations and
corresponds to a tropical $2$-plane $L$ and to a tree $T$. We define
a real valued function $\tau^d(p)$ on $\binom{[n]}{d}$ by
$\tau^d(p)_I=\sum_{\substack{ i,j \in I \\ i<j}} p_{ij}$.

\begin{prop}
$\tau^d(p)$ obeys the tropical Pl\"ucker relations.
\end{prop}

This can be checked in a routine manner, but we give a different proof
that explains the motivation behind the formula. We will call the tropical
linear space associated to $\tau^d(p)$ the $d^{\textrm{th}}$ tree space of
$L$ and denote it $\tau^d(L)$.

\begin{proof}
We know that $p_{ij}$ is realizable, meaning that we can find power series
$x_1(t)$,
\ldots $x_n(t)$, $y_1(t)$, \ldots, $y_n(t)$ in $K$ such that $p_{ij}=v(x_i y_j-x_j
y_i)$. Consider the linear space
$$\RowSpan \begin{pmatrix}
x_{1}^{d-1} & x_2^{d-1} & x_3^{d-1} & \cdots & x_n^{d-1} \\
x_{1}^{d-2} y_1 & x_2^{d-2} y_2 & x_3^{d-2} y_3 & \cdots & x_n^{d-2} y_n \\
\vdots & \vdots & \vdots & \cdots & \vdots \\
y_{1}^{d-1} & y_2^{d-1} & y_3^{d-1} & \cdots & y_n^{d-1}
\end{pmatrix}$$

The maximal minor $P_{ i_1, \ldots, i_d }$ of this matrix is the
Vandermonde determinant $\det ( x_{i_r}^{d-s}
y_{i_r}^{s-1})=\prod_{i_r<i_s} (x_{i_r} y_{i_s} - x_{i_s} y_{i_r})$. So
$$v(P_I)=\sum_{\substack{i,j \in I \\ i<j}} v(x_{i_r} y_{i_s} - x_{i_s}
y_{i_r}) =
\sum_{\substack{i,j \in I \\ i<j}} p_{ij}=\tau^d(p)_I.$$
So the $\tau^d(p)_I$ come from an actual linear space over a power series
field and hence obey the tropical Pl\"ucker relations.
\end{proof}

\begin{remark}
This construction is reminiscient of the construction of
cyclic polytopes (see, for example, \cite{Zieg} example 0.6) and
osculating flags (see \cite{Sott}, section 5), two other objects which
realize maximal combinatorics.
\end{remark}

In this section, we will give a  complete description of the bounded
part
of $\tau^d(L)$ in terms of the combinatorics of $T$. We will make heavy use
of the results of section \ref{sp}. From now on, we assume that $T$ is
trivalent.

\begin{thrm}
The vertices of $\tau^d(L)$ are of the form $\mu(T,c)^{\vee}$ where
$c$ ranges over the $\binom{n-2}{d-1}$ ways to color the internal
vertices of $T$ black and white with $d-1$ colored black and $n-d-1$
colored white. The bounded $i$-dimensional faces of $\tau^d(T)$ are in
bijection with the ordered pairs $(F,c)$, where $F$ is a forest with
$i$ trees that can be obtained by splitting $T$ along internal edges
and $c$ is a black and white coloring of the internal vertices of $F$
using $d-i$ black vertices and $n-d-i$ white vertices. $(F,c)$ is
contained in $(F',c')$ if and only if $(F',c')$ can be obtained from
$(F,c)$ by repeated splitting along edges connecting vertices of
opposite colors.
\end{thrm}

\begin{proof}
It is enough to prove the first claim, the rest then follows from the
description of the faces of $\Pi(T,c)$ in lemma \ref{faces-of-Pi(T)}.
Let $\DD$ be the subdivision of $\Delta(d,n)$
corresponding to $\tau^d(L)$, we must describe the facets of $\DD$. If $e$ is
an edge of $T$, write $[n]=A_e \sqcup B_e$ for the partition of the leaves
of $T$ induced by splitting along $e$. If we have fixed a coloring of $T$, let
$a_e$ be the number of black vertices on the $A_e$ side of $e$.

We have
$$p_I=\sum_{\substack{i,j \in I \\ i<j}} p_{ij}=(-1/2) \sum_{\substack{i,j
\in I \\ i<j}} \sum_{e \textrm{\ separates $i$\ and\ $j$}} \ell(e) =
(-1/2) \sum_e \ell(e) |A_e \cap I| \left( d- |A_e \cap I| \right).$$

For any edge $e$ of $T$, set $f_e(I)=- (1/2) \ell(e) |A_e \cap I| \left(
d- | A_e \cap I | \right)$. For $e$ containing a leaf, $f_e$ is a linear
functional. For $e$ internal, $f_e$ is convex (this uses $\ell(e)>0$). So
$\DD$ is the common refinement of the subdivisions of $\Delta(d,n)$
induced by each of the subdivisions induced by the convex functions $f_e$.
The subdivision induced by $f_e$ is to cut $\Delta(d,n)$ into the pieces
$k < \sum_{i \in A_e} x_i < k+1$ for $k \in \ZZ$.

So the facets of $\DD$ are the nonempty sets of the form
$$\left\{ x: k_e < \sum_{i \in A_e} x_i < k_e+1 \ \forall_{e \in T}
\right\}$$
for some integers $k_e$. When the
$k_e$ arise as the $a_e$ for some coloring $c$, this set is the interior
of $\Pi(T,c)$ and hence not empty. We now
must show that all the facets of $\DD$ arise in this manner. It is enough
to show that a generic point of $\Delta(d,n)$ lies in $\Pi(T,c)$ for some
coloring $c$ of $T$.

Let $x_i$ be a generic point of $\Delta(d,n)$ and let $k_e$ and $l_e$ be
the integers such that $k_e < \sum_{i \in A_e} x_i < k_e+1$, $l_e <
\sum_{i\in B_e} x_i <l_e+1$, so $k_e+l_e=d-1$. Let $v$ be an internal
vertex of $T$ with edges $e_1$, $e_2$ and $e_3$ ending at $v$. Without
loss of generality, suppose that $A_{e_s}$ is the non-$u$ side of $e_s$
for $s=1$, $2$, $3$. Then
$$k_{e_1}+k_{e_2}+k_{e_3} < \sum_{s=1}^3 \sum_{i \in A_{e_s}} x_i = d <
k_1+k_2+k_3+3$$
so $k_{e_1}+k_{e_2}+k_{e_3} = d-1$ or $d-2$. We color $v$ white if this
sum is $d-1$ and black if it is $d-2$.

We claim that precisely $d-1$ vertices are colored black. Let $B$
denote the set of vertices that are colored black and $W$ those that
are colored white. We have
$$\sum_{e \in T} (k_e+l_e) = (2n-3)(d-1)=|B| (d-2)+|W|(d-1)+n(d-1)$$
where the second expression comes from grouping the sum on
vertices. We have $|B|+|W|=n-2$ so
$$|B| (d-2)+|W|(d-1)+n(d-1)=(d-1)(2n-2)-|B|.$$
We deduce that $(2n-3)(d-1)=(2n-2)(d-1)-|B|$ and thus $|B|=d-1$. We now
will show
that $(x_i)$ is contained in the polytope $\Pi(T,c)$.

Let $e$ be an internal edge of $T$ and let $a_e$ be the number of black
vertices on a chosen side of $e$. Then $k_e < \sum_{i \in A_e} x_i <k_e
+1$ and we want to show  $a_e < \sum_{i \in A_e} x_i <a_e +1$. Let $S$ be
the subtree of $T$ lying on the $A_e$ side of $e$ and let $S$ have $s$
vertices, $s-|A_e|$ internal vertices and $s-1$ edges. We have
$$\sum_{e \subset S} (k_e+l_e)=(s-1)(d-1)=a_e (d-2)+(s-|A_e|-a_e)(d-1) +
|A_e|(d-1) - l_e$$
by once again grouping the sum on edges and on vertices. Canceling, we get
$a_e=d-1-l_e=k_e$.

So we have shown that a generic point of $\Delta(d,n)$ lies in $\Pi(T,c)$
for
some coloring $c$ of $T$ and we are done.

\end{proof}

\begin{cor}
$\tau^d(L)$ has $f_{i,d,n}=\binom{n-2i}{d-i} \binom{n-i-1}{i-1}$ bounded
faces of dimension $i$.
\end{cor}

\begin{proof}
We must count the $i$ tree forests $F$ which can be obtained by
repeated splittings of $T$ and then multiply this number by
$\binom{n-2i}{d-i}$, the number of ways to choose which vertices to
color black. Thus, it is enough to show that every trivalent tree $T$
can be split into $\binom{n-i-1}{i-1}$ different $i$ tree forests by
splitting along internal edges. Let
$F_{i,n}$ denote the number of ways to split a trivalent $n$ leaf tree
into an
$i$ tree forest. We will prove $F_{i,n}=\binom{n-i-1}{i-1}$ by
induction on $n$, it is clearly true for $n=3$.

Let $a$ and $b$ be two leaves of $T$ that have a common neighbor $v$. Let
$F$ be an $i$ tree forest obtained by splitting $T$ along internal edges.
Clearly, no sequence
of splittings can ever separate $a$ and $b$. There are then two cases: the
tree of $F$ that contains $\{ a,b \}$ contains no other leaves, or it
has some other leaf.

Case I: $\{ a,b \}$ is a component of $F$. Let $e$ be the edge joining $v$
to the rest of $T$. Then, at some point in the splitting procedure, $e$ is
split and we may as well assume that it is the first step. Let $T'$ be the
component of $T \remove e$ other than $\{ a,b \}$. The number of $F$ for
which this case applies is then the number of $i-1$ tree forests
obtainable by splitting $T'$, or $F_{i-1, n-2}$.

Case II: The component of $F$ containing $\{ a, b \}$ has additional
vertices. Let $T''$ be the tree obtained by shrinking $a$, $b$ and $v$
down to a single leaf. The number of $F$ for which this case applies is
the same as the number of $i$ tree forests obtainable by splitting $T''$,
or $F_{i,n-1}$.

So
$F_{i,n}=F_{i-1,n-2}+F_{i,n-1}=\binom{n-i-2}{i-2}+\binom{n-i-2}{i-1}=\binom{n-i-1}{i-1}$
as desired.
\end{proof}

\pagebreak
\raggedright
\pagestyle{empty}


\thebibliography{99}


\bibitem{ArdKliv}
F. Ardila and C. Klivans, ``The Bergman complex of a matroid and
phylogenetic trees'' \texttt{arxiv:math.CO/0311370}

\bibitem{Bryl}
T. Brylawski, ``A Combinatorial Model for Series-Parallel Networks''
\emph{Transactions of the AMS} \textbf{154} (1971) 1--22

\bibitem{Crapo}
H. Crapo, ``A Higher Invariant for Matroids'' \emph{Journal of
Combinatorial Theory} \textbf{2} (1967) 406--417

\bibitem{DevStur}
M. Develin and B. Sturmfels, ``Tropical Convexity'' \emph{Documenta Math.}
\textbf{9} (2004), 1--27

\bibitem{Dress1}
A. Dress and W. Wenzel ``Valuated matroids''  \emph{Advances in
Mathematics}  \textbf{93}  (1992),  no. 2, 214--250

\bibitem{Dress2}
A. Dress, ``Duality Theory for Finite and Infinite Matroids with
Coefficients''  \emph{Advances in  Mathematics} \textbf{59} 97--123 (1986)

\bibitem{Dress3}
A. Dress and W. Wenzel, ``Grassmann-Pl\"ucker Relations and Matroids with
Coefficients''  \emph{Advances in  Mathematics} \textbf{86} 68--110 (1991)

\bibitem{TreeOfLife}
A. Dress and W. Terhale, ``The Tree of Life and Other Affine Buildings'',
\emph{Documenta Mathematica}  1998,  Extra Vol. III, 565--574 (electronic)


\bibitem{GGMS}
I. Gelfand, I. Goresky, R. MacPherson, V. Serganova, ``Combinatorial
geometries, convex polyhedra, and Schubert cells'', Adv. Math. 63 (1987)
301-316

\bibitem{FirstSteps}
J. Richter-Gebert, B. Sturmfels, T. Theobald ``First Steps in Tropical Geometry'',  to appear in \emph{Idempotent Mathematics and Mathematical Physics} ed. G. Litvinov and V. Maslov, in the series Contemporary Mathematics published by the AMS. \texttt{arXiv:math.AG/0306366}


\bibitem{Sott}
F. Sottile,  ``Enumerative real algebraic geometry.  Algorithmic and
quantitative real algebraic geometry'' 139--179, DIMACS Ser. Discrete
Math. Theoret. Comput. Sci., 60, Amer. Math. Soc., Providence, RI, 2003

\bibitem{Horn}
``Horn's Problem, Vinnikov Curves and the Hive Cone'', \emph{Duke
Mathematical Journal}, to appear. \texttt{arXiv:math.AG/0311428}

\bibitem{SpeyLior}
L. Pachter and D. Speyer, ``Reconstructing Trees from Subtree Weights''
\emph{Applied Mathematics Letters} \textbf{17} (2004), p. 615 - 621.


\bibitem{TropMath}
D. Speyer and B. Sturmfels, ``Tropical Mathematics'',
\texttt{arXiv:math.CO/0408099}

\bibitem{SpeySturm}
D. Speyer and B. Sturmfels, ``The Tropical Grassmannian'' \emph{Advances
in Geometry},  Volume 4, Issue 3 (2004) p. 389 - 411


\bibitem{White1}
N. White, \emph{Theory of Matroids -- Encyclopedia of Mathematics and its
Applications vol. 26} Cambridge University Press, London (1986)


\bibitem{White2}
N. White, \emph{Matroid Applications -- Encyclopedia of Mathematics and
its Applications vol. 40} Cambridge University Press, London (1992)

\bibitem{Zieg}
G. Ziegler, \emph{Lectures on Polytopes} Graduate Texts in Mathematics
152, Springer-Verlag New York 1995
\end{document}